\documentclass[MSNbibl,number,citesort,seceqn,dvips]{arxbj}
\usepackage{upgreek}
\usepackage{graphicx}

\aid{0}
\volume{19}
\issue{5B}
\pubyear{2013}
\firstpage{2715}
\lastpage{2749}
\doi{10.3150/12-BEJ472} 

\makeatletter

\newcommand{\rref}[1]{\mbox{{\fontsize{7.6pt}{10pt}\selectfont{\ref{#1}}}}}

\renewcommand{\i}{\mathrm{i}}

\renewcommand{\d}{\mathrm{d}}

\renewcommand{\o}{\mathrm{o}}

\renewcommand{\O}{\mathrm{O}}

\newcommand{\xrightarrow}[1]{\stackrel{#1}{\hbox to 1cm{\rightarrowfill}}}
\newcommand{\xxrightarrow}[1]{\stackrel{#1}{\longrightarrow}}

\newcommand{\underset}[2]{\mathop{#2}_{#1}}

\renewcommand{\epsilon}{\varepsilon}

\newtheorem{satz}{Theorem}[section]

\newremark{rem}[satz]{Remark}

\newproclaim{assumption}[satz]{Assumption}

\newcommand{\R}{\mathbb{R}} 
\newcommand{\N}{\mathbb{N}} 
\newcommand{\E}{\mathbb{E}}

\newcommand{\cum}{\operatorname{cum}}
\newcommand{\Cov}{\operatorname{Cov}}

\newcommand{\ze}{\mathbb{Z}}

\makeatother

\begin{document}
\begin{frontmatter}

\title{A test for stationarity based on empirical processes}
\runtitle{A test for stationarity based on empirical processes}

\begin{aug}
\author{\fnms{Philip} \snm{Preu\ss}\corref{}\thanksref{e1}\ead[label=e1,mark]{philip.preuss@ruhr-uni-bochum.de}},
\author{\fnms{Mathias} \snm{Vetter}\thanksref{e2}\ead[label=e2,mark]{mathias.vetter@ruhr-uni-bochum.de}} \and
\author{\fnms{Holger} \snm{Dette}\thanksref{e3}\ead[label=e3,mark]{holger.dette@ruhr-uni-bochum.de}}
\runauthor{P. Preuÿ, M. Vetter and H. Dette} 
\address{Fakult\"at f\"ur Mathematik, Ruhr-Universit\"at Bochum,
44780 Bochum, Germany.\\ \printead{e1,e2};\\ \printead*{e3}}
\end{aug}

\pdfauthor{Philip Preu\ss, Mathias Vetter, Holger Dette}

\received{\smonth{9} \syear{2011}}
\revised{\smonth{3} \syear{2012}}

%
\begin{abstract}
In this paper we investigate the problem of testing the assumption of
stationarity in locally stationary processes. The test is based on an
estimate of a Kolmogorov--Smirnov type distance between the true time
varying spectral density and its best approximation through a
stationary spectral density. Convergence of a time varying empirical
spectral process indexed by a class of certain functions is proved, and
furthermore the consistency of a bootstrap procedure is shown which is
used to approximate the limiting distribution of the test statistic.
Compared to other methods proposed in the literature for the problem of
testing for stationarity the new approach has at least two advantages:
On one hand, the test can detect local alternatives converging to the
null hypothesis at any rate $g_T \to0$ such that $g_T T^{1/2} \to
\infty$, where $T$ denotes the sample size. On the other hand, the
estimator is based on only one regularization parameter while most
alternative procedures require two. Finite sample properties of the
method are investigated by means of a simulation study, and a
comparison with several other tests is provided which have been
proposed in the literature.
\end{abstract}

%
\begin{keyword}
\kwd{bootstrap}
\kwd{empirical spectral measure}
\kwd{goodness-of-fit tests}
\kwd{integrated periodogram}
\kwd{locally stationary process}
\kwd{non-stationary processes}
\kwd{spectral density}
\end{keyword}

\end{frontmatter}

\section{Introduction}

Most literature in time series analysis assumes that the underlying
process is second-order stationary. This assumption allows for an
elegant development of powerful
statistical methodology like parameter estimation or forecasting
techniques, but is often not justified in practice. In reality, most processes
change their second-order characteristics over time and numerous models
have been proposed to address this feature. Out of the large
literature, we mention
exemplarily the early work on this subject by Priestley \cite{priestley1965}, who
considered oscillating processes. More recently,
the concept of locally stationary processes has found considerable
attention, because in contrast to
other proposals it allows for a meaningful asymptotic theory, which is
essential for statistical inference in such models. The class of
locally stationary processes was introduced by Dahlhaus \cite{dahlhaus1996} and
particular important examples are time varying ARMA models.

While many estimation techniques for locally stationary processes were
developed (see Neumann and von Sachs \cite{neumsach1997}, Dahlhaus,
Neumann and von Sachs \cite{dahlneusachs1999}, Chiann and Morettin
\cite{changmorettin1999}, Dahlhaus and Polonik \cite{dahlpolo2006},
Dahlhaus and Subba~Rao \cite{dahlrao2006}, Van~Bellegem and von Sachs
\cite{bellegemsachs2008} or Palma and Olea \cite{palole2010} among
others), goodness-of-fit testing has found much less attention although
its importance was pointed out by many authors. von Sachs and Neumann
\cite{sacneug2000} proposed a method to test the assumption of
stationarity, which is based on the estimation of wavelet coefficients
by a localised version of the periodogram. Paparoditis
\cite{paparoditis2009} and Paparoditis \cite{paparoditis2010} used an
$L_2$ distance between the true spectral density and its best
approximation through a stationary spectral density to measure
deviations from stationarity, and most recently Dwivedi and Subba~Rao
\cite{rao2010} developed a Portmanteau type test statistic to detect
non-stationarity. However, besides the choice of a window width for the
localised periodogram which is inherent in essentially any statistical
inference for locally stationary processes, all these concepts require
the choice of at least one additional regularization parameter. For
example, the procedure proposed in Sergides and Paparoditis
\cite{sergpapa2009} relies on an additional smoothing bandwidth for the
estimation of the local spectral density. It was pointed out therein
that it is the choice of this particular tuning parameter that
influences the results of the statistical analysis substantially.

Recently, Dette, Preuss and Vetter \cite{detprevet2010} proposed a test
for stationarity which is based on an $L_2$ distance between the true
spectral density and its best stationary approximation and which does
not require the choice of that additional regularization parameter.
Roughly speaking, these authors proposed to estimate the $L_2$ distance
considered by Paparoditis \cite{paparoditis2009} by calculating
integrals of powers of the spectral density directly via Riemann sums
of the periodogram. With this idea, Dette, Preuss and Vetter
\cite{detprevet2010} avoided the integration of the smoothed
periodogram, as it was done in Paparoditis \cite{paparoditis2009} or
Paparoditis \cite{paparoditis2010}. In a comprehensive simulation study
it was shown that this method is superior compared to the other tests,
no matter how the additional smoothing bandwidths in these procedures
are chosen.

Although the test proposed by Dette, Preuss and Vetter
\cite{detprevet2010} has attractive features, it can only detect local
alternatives converging to the null hypothesis at a rate $T^{-1/4}$,
where $T$ here and throughout the paper denotes the sample size. It is
the aim of the present paper to develop a test for stationarity in
locally stationary processes which is at first able to detect
alternatives converging to the null hypothesis at the rate $g_T \to0$
such that $g_T T^{1/2} \to\infty$ and is secondly based on the concept
in Dette, Preuss and Vetter \cite{detprevet2010} for which no
additional smoothing bandwidth is needed. For this purpose, we employ a
Kolmogorov--Smirnov type test statistic to estimate a measure of
deviation from stationarity, which is defined by
\[
D:=\sup_{(v,\omega) \in[0,1]^2} {\bigl|D(v,\omega)\bigr|},
\]
where for all $(v,\omega) \in[0,1]^2$ we set
%
\begin{equation}
\label{dproc} D(v,\omega):= \frac{1}{2\uppi} \biggl( \int_0^v
\int_0^{\uppi\omega} f(u,\lambda)\,\d  \lambda \,\d u-v\int
_0^{\uppi\omega} \int_0^1
f( u,\lambda) \,\d u \,\d  \lambda \biggr)
\end{equation}
and where
$f(u,\lambda)$ denotes the time varying spectral density. Note that the
quantity $D$ is identically zero if the process is stationary, that is, if
$f(u,\lambda)$ is does not depend on $u$. The consideration of functionals
of the form (\ref{dproc}) for the construction of a test for
stationarity is natural and
was already suggested by Dahlhaus \cite{dahlhaus2009}. In particular,
Dahlhaus and Polonik \cite{dahlpolo2009} proposed an estimator
of this quantity which is based on the integrated (with respect to the
Lebesgue measure) pre-periodogram. However, in applications Riemann
sums are used
to approximate the integral and therefore the approach proposed by
these authors is not
directly implementable. In particular, it is pointed out in Example 2.7
of Dahlhaus \cite{dahlhaus2009} that
the asymptotic properties of an estimator based on Riemann
approximation have been an open problem so
far. See the discussion at the end of Section \ref{sec2} for more details.

In Section \ref{sec2}, we introduce an alternative stochastic process,
say $\{ \hat D_T (v,w)\}_{(v,w) \in[0,1]^2}$, which is based on a
summation of the localised periodogram and serves as an estimate of $\{
D (v,w)\}_{(v,w) \in[0,1]^2}$. The proposed statistic does neither
require integration of the localised periodogram with respect to an
absolutely continuous measure nor the problematic choice of a second
regularization parameter. Weak convergence of a properly standardized
version of $ \hat D_T $ to a Gaussian process is established under the
null hypothesis, local and fixed alternatives, giving a consistent
estimate of $D$. The distribution of the limiting process depends on
certain features of the data generating process which are difficult to
estimate. Therefore, the second purpose of this paper is the
development of an $\operatorname{AR}(\infty)$ bootstrap method and a proof of its
consistency. See Section \ref{sec3} for details. We also provide a
solution of the problem mentioned in the previous paragraph and prove
weak convergence of a Riemann approximation for the integrated
pre-periodogram proposed by Dahlhaus \cite{dahlhaus2009}, which is
Theorem \ref{thm1a} in the following section. As a result, we obtain
two empirical processes estimating the function $D$ defined in
(\ref{dproc}) which differ by the use of localised periodogram and
pre-periodogram in the Riemann approximations. In Section \ref{sec4},
we investigate their finite sample properties by means of a simulation
study. Although the estimator based on the pre-periodogram does not
require the specification of any regularization parameter at all, it is
demonstrated that it yields substantially less power compared to the
statistic based on the localised periodogram. Additionally, it is shown
that the latter method is extremely robust with respect to different
choices of the window width which is used for the calculation of the
localised periodogram. Moreover, we also provide a comparison with the
tests proposed in Paparoditis \cite{paparoditis2010}, Dwivedi and
Subba~Rao \cite{rao2010} and Dette, Preuss and Vetter
\cite{detprevet2010} and show that the new proposal performs better in
many situations. Finally, we present a data example, and for the sake
of a transparent presentation of the results all technical details are
deferred to the \hyperref[app]{Appendix}.

\section{The test statistic} \label{sec2}

Following Dahlhaus and Polonik \cite{dahlpolo2009}, we define a locally stationary process
via a sequence of stochastic processes $\{X_{t,T}\}_{t=1,\ldots,T}$ which
exhibit a time varying $\operatorname{MA}(\infty)$ representation, namely
%
\begin{equation}
\label{proc} X_{t,T}=\sum_{l=-\infty}^{\infty}
\psi_{t,T,l}Z_{t-l},\qquad t=1,\ldots,T,
\end{equation}
where the random variables $Z_t$ are independent identically standard
normal distributed random variables. Since the coefficients $\psi_{t,T,l}$ are in general time dependent, each process $\{ X_{t,T} \}_{t=1,\ldots,T}$ is typically not stationary. To ensure that the process
shows approximately stationary behavior on a small time interval, we
impose that there exist twice continuously differentiable functions
$\psi_l\dvtx [0,1] \rightarrow\mathbb{R}$, $l\in\mathbb{Z}$, such that
%
\begin{equation}
\label{apprbed} \sum_{l=-\infty}^{\infty}
\sup_{t = 1, \ldots,
T}\bigl|\psi_{t,T,l}-\psi_l(t/T)\bigr|=\O(1/T)
\end{equation}
as $T\to\infty$.
Furthermore, we assume that the following technical conditions
%
\begin{eqnarray}
\label{2.1a}
\sum_{l=-\infty}^{\infty} \sup_{u \in[0,1]} \bigl|
\psi_l (u)\bigr||l| &<& \infty,
\\
\label{2.1b}
\sum_{l=-\infty}^{\infty} \sup_{u \in[0,1]} \bigl|
\psi_l' (u)\bigr| &<& \infty,
\\
\label{2.1c}
\sum_{l=-\infty}^{\infty} \sup_{u \in[0,1]} \bigl|
\psi_l'' (u)\bigr| &<& \infty
\end{eqnarray}
are satisfied, which are in general rather mild. See
Dette, Preuss and Vetter \cite{detprevet2010} for a discussion. Note that variables $Z_t$ with time
varying variance $\sigma^2(t/T)$ can be included in the model by
choosing the coefficients $\psi_{t,T,l}$ in (\ref{proc}) appropriately.

Set
\[
\psi\bigl(u,\exp(-\i  \lambda)\bigr):= \sum_{l=-\infty}^{\infty}
\psi_l(u) \exp(-\i  \lambda l).
\]
Then the function
\[
f(u,\lambda)=\frac{1}{2 \uppi}\bigl|\psi\bigl(u,\exp(-\i  \lambda)\bigr)\bigr|^2
\]
is well defined and called the time varying spectral density of $\{
{X}_{t,T} \}_{t=1,\ldots,T}$, see Dahlhaus \cite{dahlhaus1996}. It is
continuous by assumption and can roughly be estimated by a local
periodogram. To be precise, we assume without loss of generality that
the total sample size $T$ can be decomposed as $T=NM$, where $N$ and
$M$ are integers and $N$ is even. Furthermore, we define
\[
I_N^X(u,\lambda):=\frac{1}{2 \uppi N} \Biggl|\sum
_{s=0}^{N-1} X_{\lfloor
uT \rfloor-N/2+1+s,T} \exp(-\i  \lambda s)
\Biggr|^2,
\]
which is the local periodogram at time $u$ proposed by
Dahlhaus \cite{dahlhaus1997}. Here, we have set
$X_{j,T} =0$, if $j \notin\{ 1,\ldots, T\}$. This is the usual
periodogram computed from the observations $X_{\lfloor uT
\rfloor-N/2+1,T}, \ldots, X_{\lfloor uT \rfloor+N/2,T}$. The arguments
employed in the \hyperref[app]{Appendix} show that
\[
\E\bigl(I_N^X(u,\lambda)\bigr)=f(u,
\lambda)+\O(1/N)+\O(N/T),
\]
and therefore the statistic $I_N^X(u,\lambda)$ is an asymptotically
unbiased estimator for the spectral density if $N \rightarrow\infty$
and $N=\o(T)$. However, $I_N^X(u,\lambda)$ is not consistent just as the
usual periodogram.

We now consider an empirical version of the function $D(v,\omega)$
defined in (\ref{dproc}), that is,
%
\begin{equation}
\label{hatD} \hat D_T(v,\omega):=\frac{1}{T}\sum
_{j=1}^{\lfloor vM \rfloor}\sum_{k=1}^{\lfloor\omega{N}/{2} \rfloor}I_N^X(u_j,
\lambda_k)-\frac
{\lfloor vM \rfloor}{M}\frac{1}{T}\sum
_{j=1}^M \sum_{k=1}^{\lfloor
\omega{N}/{2} \rfloor}I_N^X(u_j,
\lambda_k),
\end{equation}
where the points
\[
u_j:=\frac{t_j}{T}:=\frac{N(j-1)+N/2}{T},\qquad j=1,\ldots,M,
\]
define an equidistant grid of the interval $[0,1]$ and
\[
\lambda_k:=\frac{2\uppi k}{N},\qquad k=1,\ldots,\frac{N}{2},
\]
denote the Fourier frequencies.
It follows from the proof of Theorem \ref{thm1} in the \hyperref[app]{Appendix} that
for every $v \in[0,1]$ and $\omega\in[0,1]$ we have
\begin{eqnarray*}
\E\bigl(\hat D_T(v,\omega)\bigr)
&=& \frac{1}{T}\sum
_{j=1}^{\lfloor vM \rfloor
}\sum_{k=1}^{\lfloor\omega{N}/{2} \rfloor}f(u_j,
\lambda_k)\\
&&{}-\frac
{\lfloor vM \rfloor}{M}\frac{1}{T}\sum
_{j=1}^M \sum_{k=1}^{\lfloor
\omega{N}/{2} \rfloor}f(u_j,
\lambda_k)+\O(1/N)+\O\bigl(N^2/T^2\bigr)
\\
&=& D(v,\omega)+\O(1/N)+\O(N/T),
\end{eqnarray*}
where the latter identity is due to the approximation error of the
Riemann sum. This error can be improved, if we replace $D(v,\omega)$ by
its discrete time approximation, that is,
\[
D_{N,M}(v,\omega):=D \biggl( \frac{\lfloor vM \rfloor}{M},\frac
{\lfloor
\omega{N}/{2} \rfloor}{{N}/{2}}
\biggr)
\]
for which the representation
%
\begin{equation}
\label{expan} \E\bigl(\hat D_T(v,\omega)\bigr)=D_{N,M}(v,
\omega)+\O(1/N)+\O\bigl(N^2/T^2\bigr)
\end{equation}
holds. The approximation error of the Riemann sum in (\ref{expan})
becomes smaller due to the choice of the midpoints $u_j$. The rate of
convergence will be $T^{-1/2}$ later on, so we need the $\O(\cdot
)$-terms to vanish asymptotically after multiplication with $\sqrt{T}$.
Therefore, we define an empirical spectral process by
\[
\hat G_T(v,\omega):=\sqrt{T} \Biggl( \frac{1}{T}\sum
_{j=1}^{\lfloor vM
\rfloor}\sum_{k=1}^{\lfloor\omega{N}/{2} \rfloor
}I_N^X(u_j,
\lambda_k)-\frac{\lfloor vM \rfloor}{M}\frac{1}{T}\sum
_{j=1}^M \sum_{k=1}^{\lfloor\omega{N}/{2} \rfloor}I_N^X(u_j,
\lambda_k) -D_{N,M}(v,\omega) \Biggr)
\]
and assume
%
\begin{equation}
\label{assNM} N \rightarrow\infty,\qquad M \rightarrow\infty,\qquad \frac
{T^{1/2}}{N}
\rightarrow0,\qquad \frac{N}{T^{3/4}} \rightarrow0.
\end{equation}
Our first result specifies the asymptotic properties of the empirical
process\break  $(\hat G_T(v,\omega))_{(v,\omega)\in[0,1]^2}$, both
under the null hypothesis and under a fixed alternative. The null
hypothesis of stationarity is formulated as
%
\begin{equation}
\label{H0} H_0\dvtx f(u,\lambda) \mbox{ is independent of } u,
\end{equation}
which is a little different from genuine second-order stationarity,
since it only means that the coefficients $\psi_{t,T,l}$ in (\ref
{proc}) can be approximated by time independent terms $\psi_l$. Thanks
to the continuity of the time varying spectral density, the alternative
corresponds to the property that there is some $\lambda$ such that $u
\mapsto f(u, \lambda)$ is not a constant function. Finally, the symbol
$\Rightarrow$ denotes weak
convergence in $[0,1]^2$.
%
\begin{satz}
\label{thm1} Suppose we have a locally stationary process as defined in
(\ref{proc}) with independent and standard normal innovations $Z_t$.
Furthermore assume that the assumptions (\ref{apprbed})--(\ref{2.1c}) and
(\ref{assNM}) are satisfied. Then as $T \rightarrow\infty$ we have
%
\begin{equation}
\label{weak} \bigl(\hat G_T(v,\omega)\bigr)_{(v, \omega) \in[0,1]^2}
\Rightarrow \bigl(G(v,\omega )\bigr)_{(v, \omega) \in[0,1]^2},
\end{equation}
where $ (G(v,\omega))_{(v, \omega) \in[0,1]^2}$
is a Gaussian process with mean zero and covariance structure
\begin{eqnarray*}
&&\Cov\bigl(G(v_1,\omega_1),G(v_2,
\omega_2)\bigr)
\\
&&\quad=\frac{1}{2\uppi}\int_0^1\int
_0^{\uppi\min(\omega_1,\omega_2)} \bigl(1_{[0,v_1]}(u)-v_1
\bigr) \bigl(1_{[0,v_2]}(u)-v_2\bigr)f^2(u,\lambda) \,\d
\lambda \,\d u.
\end{eqnarray*}
\end{satz}

Under the null hypothesis, we have $D_{N,M}(v,\omega)=0$ for all $N,M
\in\N$ and for all $v,\omega\in[0,1]$. Therefore, we obtain
\[
\bigl(\sqrt{T} \hat D_T(v,\omega)\bigr)_{(v,\omega) \in[0,1]^2} \Rightarrow
\bigl(G(v,\omega)\bigr)_{(v,\omega) \in[0,1]^2},
\]
which yields
%
\begin{equation}
\label{dnull} \sqrt{T} \sup_{(v,\omega) \in[0,1]^2} \bigl|\hat D_T(v,\omega)\bigr|
\xxrightarrow {D} \sup_{(v,\omega) \in[0,1]^2} \bigl|G(v,\omega)\bigr|
\end{equation}
under the null hypothesis (\ref{H0}).
An asymptotic level $\alpha$ test is then obtained by rejecting the
null hypothesis of stationarity
whenever $\sqrt{T} \sup_{(v, \omega) \in[0,1]^2} |\hat D_T(v,\omega)|$
exceeds the ($1-\alpha$)\% quantile of
the distribution of the random variable $\sup_{(v,\omega)\in[0,1]^2}
|G(v,\omega)|$.
On the other hand, under the alternative there is a pair $(v, \omega)$
such that $D(v,\omega) \neq0$. The fact that $D_{N,M}$ converges
uniformly to $D$ together with Theorem \ref{thm1} yields consistency of
this test.
Note also that under the null hypothesis $H_0$
the covariance structure of the limiting process in Theorem \ref{thm1}
simplifies to
%
\begin{equation}
\label{covnull} \Cov\bigl(G(v_1,\omega_1),G(v_2,
\omega_2)\bigr)=\frac{\min( v_1,v_2) -
v_1v_2}{2\uppi} \int_0^{\uppi\min(\omega_1,\omega_2)}f^2(
\lambda) \,\d \lambda
\end{equation}
and depends on the unknown spectral density $f$. In order to avoid the
estimation of the integral over the squared spectral density, we
propose to approximate the quantiles of the limiting distribution by an
$\operatorname{AR}(\infty)$ bootstrap,
which will be described in the following section.

An alternative estimator for the time varying spectral
density is given by
\[
J_T(u,\lambda):=\frac{1}{2\uppi} \sum_{k\dvtx 1 \leq\lfloor uT+1/2\pm k/2
\rfloor\leq T}
X_{\lfloor uT+1/2+k/2 \rfloor}X_{\lfloor uT+1/2-k/2
\rfloor} \exp(-\i  \lambda k),
\]
which is called the pre-periodogram (see Neumann and von Sachs \cite{neumsach1997}). As for
the usual periodogram, it is asymptotically unbiased, but again not
consistent. Based on this statistic, we define an alternative
process by
%
\begin{eqnarray}
\label{preperio1}
\hat H_T^1(v,\omega)&:=&\sqrt{T} \Biggl(
\frac{1}{T^2}\sum_{j=1}^{\lfloor
vT \rfloor}\sum
_{k=1}^{\lfloor\omega{T}/{2} \rfloor}J_T(j/T,
\lambda_{k,T})\nonumber\\[-8pt]\\[-8pt]
&&\hspace*{22pt}{}- \frac{\lfloor vT \rfloor}{T^3}\sum_{j=1}^T
\sum_{k=1}^{\lfloor
\omega
{T}/{2} \rfloor}J_T(j/T,
\lambda_{k,T}) -D(v,\omega) \Biggr),\nonumber
\end{eqnarray}
where the Fourier frequencies become $\lambda_{k,T} = {2\uppi k}/{T}$
now. Convergence of the finite dimensional distributions of the process
$(H_T^1(v,\omega))_{(v, \omega) \in[0,1]^2}$ to the ones of the
limiting process $(G(v,\omega))_{(v,\omega) \in[0,1]^2}$ has already
been shown in Dahlhaus \cite{dahlhaus2009}. Tightness can be shown
using similar arguments as given in the \hyperref[app]{Appendix} for
the proof of Theorem \ref{thm1}, which are not stated here for the sake
of brevity. As a consequence, we obtain the following result.
%
\begin{satz}
\label{thm1a} If the assumptions of Theorem \ref{thm1}
are satisfied, then as $T \rightarrow\infty$ we have
\[
\bigl(\hat H_T^1(v,\omega)\bigr)_{(v, \omega) \in[0,1]^2}
\Rightarrow \bigl(G(v,\omega)\bigr)_{(v, \omega) \in[0,1]^2},
\]
where $ (G(v,\omega))_{(v, \omega) \in[0,1]^2}$
is the Gaussian process defined in Theorem \ref{thm1}.
\end{satz}

Because the use of $\hat H_T^1(v,\omega)$ instead of $\hat
G_T(v,\omega
)$ does not require the choice of the quantity $N$, which specifies the
number of observations used
for the calculation of the local periodogram, it might be appealing to
construct a
Kolmogorov--Smirnov type test for stationarity
on the basis of this process.
However, we will demonstrate in Section \ref{sec4} by means of a
simulation study that
for realistic sample sizes
the method which employs the pre-periodogram
is clearly outperformed by the approach based on the local periodogram.
Our numerical results also show that the use of the local
periodogram\vadjust{\goodbreak}
is not very sensitive with respect
to the choice
of the regularization parameter $N$ either, and therefore we strictly
recommend to use the
latter approach when constructing a Kolmogorov--Smirnov test.
%
\begin{rem}
\label{rem1}
The convergence of a modified version of the process
(\ref{preperio1}) to the limiting Gaussian process $(G(v,\omega))_{(v,
\omega) \in[0,1]^2}$ of Theorem \ref{thm1} was shown in
Dahlhaus and Polonik \cite{dahlpolo2009},
where the Riemann sum over the Fourier frequencies was replaced by the
integral with respect to the Lebesgue measure. More precisely,
these authors considered the process
\begin{eqnarray*}
\bigl(\hat H_T^2(v,\omega)\bigr)_{(v, \omega) \in[0,1]^2}&:=&
\frac{1}{2\uppi
\sqrt {T} } \Biggl( \sum_{j=1}^{\lfloor vT \rfloor}
\int_0^{\uppi\omega} J_T(j/T,\lambda)\,\d  \lambda
\\
&&\hspace*{35.3pt}{} - {v}\sum_{j=1}^T
\int_0^{\uppi\omega} J_T(j/T, \lambda) \,\d
\lambda -D(v,\omega) \Biggr)_{(v, \omega) \in[0,1]^2}
\end{eqnarray*}
instead of $(H_T^1(v,\omega))_{(v, \omega) \in[0,1]^2}$ and proved its
weak convergence. This is a rather typical result, as many other
asymptotic results are only shown for the integral (instead of the sum
over the Fourier coefficients) over the local periodogram or the
pre-periodogram; see, for example, Dahlhaus \cite{dahlhaus1997} or
Paparoditis \cite{paparoditis2010}. The transition from these results to analogue
statements for the corresponding Riemann approximations is by no means
obvious. For example, although it is appealing to assume that
\[
\int_0^\uppi I_N^X(u,
\lambda) \,\d \lambda= \frac{2\uppi}{N} \sum_{k=1}^{
{N}/{2}}
I_N^X(u,\lambda_k)+\O(1/N)
\]
holds because of the Riemann approximation error, this identity is in
general not valid, as the derivative ${\partial I_N^X(u,\lambda
)}/{\partial\lambda}$ is not uniformly bounded in $N$. A demonstrative
explanation of this fact is that $I_N^X(u,\lambda_{k_1})$ and
$I_N^X(u,\lambda_{k_2})$ are asymptotically independent whenever
$k_1\not= k_2$. Thus in general asymptotic results for integrated local
periodogram or pre-periodogram cannot be directly transferred to the
corresponding Riemann approximations.
These difficulties were also explicitly pointed out in Example 2.7 of
Dahlhaus \cite{dahlhaus2009}. Note further that asymptotic tightness has neither
been studied for an integrated nor for a summarized local periodogram
in the literature so far.
\end{rem}
%
\begin{rem}
\label{rem2}
Suppose that we are in the situation of local alternatives, that is, we have
%
\begin{equation}
\label{loc} f_T(u,\lambda)=f(\lambda)+g_T k (u,\lambda)
\end{equation}
for some deterministic sequence $g_T$ and an appropriate function $k$
such that (\ref{loc}) defines a time varying spectral density. Note
that a locally stationary process with this specific spectral density
can easily be constructed through the equation
\[
X_{t,T}=\int_{-\uppi}^{\uppi} \exp(\i \lambda
t)A_T(t/T,\lambda)\,\d \xi (\lambda),
\]
where $\xi$ is an orthogonal increment Gaussian process and
$A_T(u,\lambda)$ is a function such that $f_T(u,\lambda
)=|A_T(u,\lambda
)|^2$. See Dahlhaus \cite{dahlhaus1997}.

A careful inspection of the proofs in the \hyperref[app]{Appendix} shows that (\ref
{weak}) with centering term $D_{N,M}(v, \omega)=0$ and asymptotic
covariance (\ref{covnull}) also holds
in the case where $g_T = \o(1/ \sqrt{T}) $.
Moreover, if $g_T={1}/{\sqrt{T}}$ an analogue of
Theorem \ref{thm1} can be obtained where
the centering term $D_{N,M} (v,\omega)$ in the definition of $\hat
G_T(v,\omega)$ is replaced
by
\begin{eqnarray*}
D_{N,M, k} (v,\omega) &=& \frac{1}{2\uppi\sqrt{T}} \biggl( \int
_0^ {\lfloor vM \rfloor /M} \int
_0^{{2 \uppi\lfloor\omega{N}/{2} \rfloor}/{{N}}
} k(u,\lambda)\,\d  \lambda \,\d u\\
&&\hspace*{35pt}{} -
\frac{\lfloor vM \rfloor}{M} \int_0^{
{2 \uppi\lfloor\omega{N/2} \rfloor}/{{N}}
} \int
_0^1 k( u,\lambda) \,\d u \,\d  \lambda \biggr),
\end{eqnarray*}
which is the original $D_{N,M}$ but with $T^{-1/2}k(u,\lambda)$ playing
the role of $f(u,\lambda)$. In this
case, the appropriately centered process converges weakly to a Gaussian
process $\{ G(v,\omega)\}_{(v, \omega)_\in[0,1]^2}$ with
covariance structure given by (\ref{covnull}) as well. A similar
comment applies to the process $\hat H_T^1$ defined in (\ref{preperio1}).
This means that the tests based on the processes $ \hat G_T$ and $\hat
H_T^1$ can detect alternatives converging to the null hypothesis
at any rate $g_T \to0$ such that $g_T T^{1/2} \to\infty$.
In contrast, the proposal of Dette, Preuss and Vetter \cite{detprevet2010} is based on an $L_2$ distance
between $f(u,\lambda)$ and $\int_0^1f(v,\lambda)\,\d v$ and is therefore
only able to detect
alternatives converging to the null hypothesis at a rate $T^{-1/4}$.
\end{rem}
%
\begin{rem}
\label{rem3}
In Theorems \ref{thm1} and \ref{thm1a}, we assume the existence
of second order derivatives for the approximating functions $\psi_l(u)$. Nevertheless, it is straightforward to show that our test also
detects fixed alternatives in which the $\psi_l(u)$ admit a finite
number of points
of discontinuity. We furthermore conjecture that the constraints in
Theorems \ref{thm1} and \ref{thm1a} can be weakened to some kind
of condition on the total variation of $\psi_l(u)$ as in Definition 2.1
in Dahlhaus and Polonik~\cite{dahlpolo2009}.
\end{rem}

\section{Bootstrapping the test statistic} \label{sec3}

To approximate the limiting distribution of $\sup_{(v,\omega) \in
[0,1]^2} |G(v,\omega)|$, we will employ an $\operatorname{AR}(\infty)$ bootstrap
approximation, which was introduced by Krei{\ss} \cite{Kreiss1988}. To ensure
consistency of the bootstrap procedure described later, we have to
consider the stationary process $Y_t$ with spectral density $\lambda
\mapsto\int_0^1 f(u,\lambda) \,\d u$ first, which coincides with $X_{t,T}$
in case the latter process is stationary. We have to impose the
following main assumption.
%
\begin{assumption}\label{ass2}
We assume that the spectral density $\lambda\mapsto\int_0^1f(u,\lambda
) \,\d u$ is strictly positive and that the process $Y_t$ has an
$\operatorname{AR}(\infty )$ representation, that is,
%
\begin{equation}
\label{statproc} Y_t=\sum_{j=1}^\infty
a_j Y_{t-j}+Z_t^{\mathrm{AR}},
\end{equation}
where $(Z_j^{\mathrm{AR}})_{j\in\ze}$ denotes a Gaussian white noise process
with some variance $\sigma^2>0$ and the sequence $(a_j)_{j \in\mathbb
{N}}$ of coefficients satisfies $\sum_{j=1}^\infty|a_j|<\infty$ and
%
\begin{equation}
\label{gl1} 1-\sum_{j=1}^\infty
a_jz^j \not= 0 \qquad\mbox{for } |z|\leq1.
\end{equation}
\end{assumption}
Note that $(Y_t)_{t \in\mathbb{Z}}$ possesses an $\operatorname{MA}(\infty)$ representation
%
\begin{equation}
\label{procstat} Y_t=\sum_{l=-\infty}^\infty
\psi_l Z_{t-l},
\end{equation}
where the $Z_t$ are the same as in (\ref{proc}) and the $\psi_l$ are
some appropriately defined constants. The random variables $Z_t^{\mathrm{AR}}$
in (\ref{statproc}) do not necessarily coincide with the $Z_t$ from
(\ref{procstat}), even though this could be ensured by assuming that
the $\operatorname{MA}(\infty)$ representation in (\ref{procstat}) corresponds to the
Wold representation of $Y_t$. See, for example, Kreiss, Paparoditis and
Politis \cite{kreisspappol2011} for a comprehensive illustration.

We have to introduce a second class of stationary processes, namely
$(Y_t^{\mathrm{AR}}(p))_{t \in\mathbb{Z}}$ for arbitrary integer $p$, which is
the process defined through
%
\begin{equation}
\label{arpstat} Y_t^{\mathrm{AR}}(p)=\sum
_{j=1}^p a_{j,p}Y_{t-j}^{\mathrm{AR}}(p)+Z_t^{\mathrm{AR}}(p),
\end{equation}
where
%
\begin{equation}
\label{ar2} (a_{1,p},\ldots,a_{p,p}):= \underset{b_{1,p},\ldots,b_{p,p}}
{\operatorname{argmin}} \E \Biggl(Y_{t}-\sum
_{j=1}^pb_{j,p}Y_{t-j}
\Biggr)^2
\end{equation}
and $(Z_t^{\mathrm{AR}}(p))_{t \in\mathbb{Z}}$ is a Gaussian white noise
process with mean zero and variance
\[
\sigma_p^2=\E \Biggl(Y_t-\sum
_{j=1}^p a_{j,p}Y_{t-j}
\Biggr)^2.
\]
In other words, $Y_t^{\mathrm{AR}}(p)$ corresponds to the best
$\operatorname{AR}(p)$ model which can be fitted to the process~$Y_t$.
Lemma 2.2 in Kreiss, Paparoditis and Politis \cite{kreisspappol2011}
ensures that for growing $p$
%
\begin{equation}\label{arbound2}
\sum_{k=1}^p (1+k) |a_{k,p}-a_k|
\rightarrow0,
\end{equation}
thus the process $Y_t^{\mathrm{AR}}(p)$ becomes `close' to the process $Y_t$.

The bootstrap procedure now works by fitting an $\operatorname{AR}(p)$ model to the
observed data $X_{1,T},\ldots,X_{T,T}$, where the parameter $p=p(T)$
increases with the sample size $T$. To be precise, we first calculate
an estimator $(\hat a_{1,p,T},\ldots,\hat a_{p,p,T})$ for
%
\begin{equation}
\label{ar} (a_{1,p,T},\ldots,a_{p,p,T})= \underset{b_{1,p,T},\ldots,b_{p,p,T}}
{\operatorname{argmin}} \E \Biggl(X_{t,T}-\sum
_{j=1}^p b_{j,p,T}X_{t-j,T}
\Biggr)^2
\end{equation}
and then simulate a pseudo series $X_{1,T}^*,\ldots,X_{T,T}^*$ according
to the model
\begin{eqnarray*}
X_{t,T}^*&=&X_{t,T};\qquad t=1,\ldots,p,
\\
X_{t,T}^*&=&\sum_{j=1}^p \hat
a_{j,p,T}X_{t-j,T}^*+Z_j^*;\qquad p<t \leq T.
\end{eqnarray*}
Here, the quantities $Z_j^*$ denote independent and normal distributed
random variables with mean zero and variance
%
\begin{equation}
\label{sigdach} \hat\sigma_p^2:=\frac{1}{T-p} \sum
_{t=p+1}^{T}(\hat z_t-\overline
z_T)^2,
\end{equation}
where $\overline z_T:=\frac{1}{T-p} \sum_{t=p+1}^T \hat z_t$
and
\[
\hat z_t:=X_{t,T}-\sum_{j=1}^p
\hat a_{j,p,T} X_{t-j,T} \qquad\mbox{for } t=p+1,\ldots,T,
\]
thus $\hat\sigma_p^2$ is the standard variance estimator of the error
process $\hat z_t$. We now define
the statistic $\hat G_T^*(v,\omega)$ in the same way as $\hat
G_T(v,\omega)$ where the original observations
$X_{1,T},\ldots,X_{T,T}$ are replaced by the bootstrap replicates
$X_{1,T}^*,\ldots,X_{T,T}^*$. To assure that this procedure approximates
the limiting distribution corresponding to the null hypothesis both
under the null hypothesis and the alternative, we need the following
technical conditions:
%
\begin{assumption}
\label{ass1}
\begin{enumerate}[(iii)]
\item[(i)] $p=p(T)\in[p_{\mathrm{min}}(T),p_{\mathrm{max}}(T)]$, where
$
p_{\mathrm{max}}(T) \geq p_{\mathrm{min}}(T) \xrightarrow{T \rightarrow\infty
} \infty
$ and
%
\begin{equation}
\label{asspNM} \frac{p_{\mathrm{max}}^3(T) \sqrt{\log(T)}}{\sqrt{T}} =\O(1).
\end{equation}
\item[(ii)] The estimators for the AR parameters defined by (\ref
{ar}) satisfy
%
\begin{equation}
\label{arestbound} \max_{1\leq j \leq p}|\hat a_{j,p,T}-a_{j,p}|=\O
\bigl(\sqrt{\log(T)/T}\bigr),
\end{equation}
uniformly with respect to $p \leq p(T)$.
\item[(iii)] The estimate $\hat\sigma_p^2$ defined in (\ref{sigdach})
converges in probability to $\sigma^2>0$.
\end{enumerate}
\end{assumption}
All assumptions are rather standard in the framework of an $\operatorname{AR}(\infty)$
bootstrap; see, for example, Krei{\ss} \cite{kreiss1997} or Berg,
Paparoditis and Politis \cite{bergpappolitis2010}. Thanks to
(\ref{asspNM}), assumption (\ref {arestbound}) is, for example,
satisfied for the least squares or the Yule--Walker estimators; see
Hannan and Kavalieris~\cite{hankav1986}. The latter condition is
extremely important, as it implies that $X_{t,T}^*$ shows a similar
behavior as the $\operatorname{AR}(p)$ process $Y_t^{\mathrm{AR}}(p)$ and is therefore also
`close' to $Y_t$ in a similar sense as~(\ref{arbound2}). Therefore, we
can expect that statistics based on the bootstrap replicates behave in
the same way as those based on a stationary process. Precisely, we
obtain the following result which implies consistency of the bootstrap
procedure described above.
%
\begin{satz}
\label{thm2} Suppose that the assumptions of Theorem \ref{thm1} hold
and that furthermore Assumptions \ref{ass2} and \ref{ass1} are
satisfied. Then as $T \rightarrow\infty$ we have conditionally on
$X_{1,T},\ldots,X_{T,T}$
\[
\bigl(\hat G_T^*(v,\omega)\bigr)_{(v, \omega) \in[0,1]^2} \Rightarrow\bigl(
\tilde G(v,\omega)\bigr)_{v \in[0,1], \omega\in[0,1]},
\]
where $(\tilde G(v,\omega))_{(v, \omega) \in[0,1]^2} $ denotes a
centered Gaussian process with covariance structure
\[
\Cov\bigl(\tilde G(v_1,\omega_1),\tilde
G(v_2,\omega_2)\bigr)=\frac{\min(
v_1,v_2) - v_1v_2}{2\uppi} \int
_0^{\uppi\min(\omega_1,\omega_2)} \biggl( \int_0^1
f(u, \lambda ) \,\d u \biggr)^2 \,\d \lambda.
\]
\end{satz}
We now obtain empirical quantiles of $\sup_{(v,\omega) \in[0,1]^2}
|G(v,\omega)|$ by calculating $\hat D_{T,i}^*:=\sup_{(v,\omega) \in
[0,1]^2} |\hat G_{T,i}^*(v,\omega)|$ for $i=1,\ldots,B$ where $\hat
G_{T,1}^*(v,\omega),\ldots,\hat G_{T,B}^*(v,\omega)$ are the $B$ bootstrap
replicates of $\hat G_T(v,\omega)$. The null hypothesis is then
rejected, whenever
%
\begin{equation}
\label{test} \sqrt{T} \sup_{(v,\omega) \in[0,1]^2} \bigl|\hat D_T(v,\omega)\bigr|>
\bigl(\hat D_T^*\bigr)_{T,\lfloor(1-\alpha)B \rfloor},
\end{equation}
where $(\hat D_T^*)_{T,1},\ldots,(\hat D_T^*)_{T,B}$ denotes the order
statistic of $\hat D_{T,1}^*,\ldots,\hat D_{T,B}^*$. The test has
asymptotic level $\alpha$ because of Theorem \ref{thm2} and is
consistent within the class of alternatives satisfying Assumptions \ref
{ass2} and \ref{ass1}. This follows, since conditionally on
$X_{1,T},\ldots,X_{T,T}$ each bootstrap statistic
$\sup_{(v,\omega) \in[0,1]^2} | \hat G_T^*(v,\omega)| $ converges to
a non-degenerate random variable,
while $\sqrt{T} \sup_{(v,\omega)\in[0,1]^2} | \hat D_T(v,\omega) |$
converges to infinity by Theorem \ref{thm1}.
We finally point out that similar results can be shown for the
statistic which is obtained by replacing
the localised periodogram in $\hat D_T$ by the pre-periodogram. The
technical details are omitted for the sake of brevity, but
the finite sample performance of this alternative approach will be
investigated in the following section as well.

\section{Finite sample properties} \label{sec4}

\subsection{Choosing the parameter}

We first comment on how to choose the parameters $N$ and $p$ in
concrete applications. Although the proposed method does not show much
sensitivity with respect to different choices of both parameters,\vadjust{\goodbreak} we
select $p$ throughout this section as the minimizer of the AIC
criterion dating back to Akaike \cite{akaike1973}, which is defined by
\[
\hat p = \mathop{\operatorname{argmin}}_p \frac{1}{T} \sum
_{k=1}^{{T}/{2}} \biggl( \log\bigl(f_{\hat\theta(p)}(
\lambda_{k,T})\bigr)+\frac{I_T^X(\lambda_{k,T})}{f_{\hat\theta(p)}(\lambda_{k,T})} \biggr) + p/T
\]
in the context of stationary processes. See also Whittle \cite{whittle1} and
Whittle \cite{whittle2}. Here, $f_{\hat\theta(p)}$ is the spectral density of
a stationary $\operatorname{AR}(p)$ process with the fitted coefficients and $I_T^X$
is the usual stationary periodogram. Therefore, we focus in the
following discussion
on the sensitivity analysis of the test (\ref{test}) with respect to
different choices of $N$, and we will see that the particular choice of
that tuning parameter has typically very little influence on the
outcome of the test.

\subsection{Bootstrap approximation}

Let us illustrate now how well the proposed bootstrap method
approximates the distribution of
the statistic $\sqrt{T} \sup_{(v,\omega)\in[0,1]^2} |\hat
D_T(v,\omega
)|$ under the null hypothesis. For this purpose, we simulate
observations from the stationary $\operatorname{AR}(1)$ model
%
\begin{equation}
\label{simproc} X_{t,T}=0.5X_{t-1,T}+Z_t,\qquad t=1,\ldots,T,
\end{equation}
for $T=128$. In particular, we generate 1000 versions of this process
and calculate each time the test statistic $\sqrt{T} \sup_{(v,\omega)
\in[0,1]^2} |\hat D_T(v,\omega)|$, both for $N=16$ and $N=8$. These
outcomes can be used to estimate the exact distribution of the test
statistic. In a next step, we choose randomly 10 series from the 1000
replications of (\ref{simproc}), for which we calculate another $1000$
%
\begin{figure}

\includegraphics{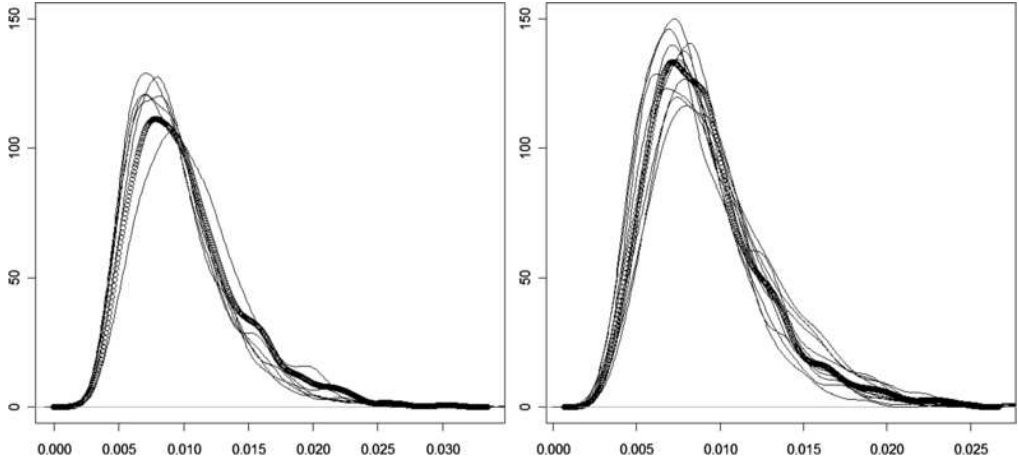}

\caption{Estimated densities\vspace*{1pt} of the distribution of the
statistic $\sqrt{T} \sup_{(v,\omega)\in[0,1]^2} |\hat D_T(v,\omega)|$
under the null hypothesis. The dotted line is the estimated exact
density while the solid lines corresponds to the estimated densities of
the bootstrap approximations. Left panel: $N=8$; right panel: $N=16$.}
\label{figure1}
\end{figure}
bootstrap approximations each. Based on these bootstrap replications,
we estimate the density of the corresponding bootstrap approximations
of the test statistic as well. The plots comparing these densities are
given in Figure \ref{figure1} where the dotted line corresponds to the
estimated exact density while the dashed lines show the $10$ estimated
densities of the bootstrap approximations.

\subsection{Size and power of the test}

In this section, we investigate the size and power of the test (\ref
{test}) and the analogue based on the pre-periodogram. We also
compare these methods with three other tests for stationarity, which
recently have been proposed in the literature.
All reported results are based
on $200$ bootstrap replications and $1000$ simulation runs under the
null hypothesis while we use $500$ simulation runs under the
alternative. To study the approximation of the nominal level, we
simulate $\operatorname{AR}(1)$ processes
%
\begin{equation}
\label{AR} X_t=\phi X_{t-1}+ Z_t,\qquad t\in\ze,
\end{equation}
and $\operatorname{MA}(1)$ processes
%
\begin{equation}
\label{MA} X_t=Z_t+\theta Z_{t-1},\qquad t\in\ze,
\end{equation}
for different values of the parameters $\phi$ and $\theta$, where the
$Z_t$ are independent and standard normal distributed random variables
throughout the whole section.
The corresponding results are depicted in Tables \ref{Table1} and \ref
{Table1a}, respectively,
and we observe a precise approximation of the nominal level in the
$\operatorname{AR}(1)$ case for $\phi\in\{-0.5,0,0.5,0.9\}$ and in the $\operatorname{MA}(1)$ case
%
\begin{table}
\caption{Rejection probabilities of the test (\protect\ref
{test}) under
the null hypothesis. The data
was generated according to model (\protect\ref{AR})}
\label{Table1}
\begin{tabular*}{\tablewidth}{@{\extracolsep{\fill}}lllllllllllll@{}}
\hline
& &  & \multicolumn{2}{l}{$\phi=-0.9$} &
\multicolumn{2}{l}{$\phi=-0.5$}
& \multicolumn{2}{l}{$\phi=0$} & \multicolumn{2}{l}{$\phi=0.5$}
& \multicolumn{2}{l@{}}{$\phi=0.9$} \\[-4pt]
& & & \multicolumn{2}{c}{\hrulefill} &
\multicolumn{2}{c}{\hrulefill} & \multicolumn{2}{c}{\hrulefill}
& \multicolumn{2}{c}{\hrulefill} & \multicolumn{2}{c@{}}{\hrulefill}\\
$T$&$N$&$M$& $5 \%$ & $10 \%$ & $5 \%$ & $10 \%$ & $5 \%$ & $10 \%$ &
$5 \%$ & $10\%$ & $5 \%$ & $10\%$ \\ \hline
\hphantom{0}$64$&\phantom{0}$8$&\phantom{0}$8$& $0.021$ & $0.069$ & $0.025$&$0.060$
&$0.035$&$0.086$&$0.050$&$0.099$ & $0.044$ & $0.108$ \\
$128$&$16$&\phantom{0}$8$& $0.022$ & $0.063$ &$0.031$&$0.077$
&$0.042$&$0.081$&$0.034$&$0.092$ & $0.050$ & $0.099$ \\
$128$&\phantom{0}$8$&$16$& $0.020$ & $0.066$ &$0.030$&$0.076$
&$0.038$&$0.083$&$0.055$&$0.102$ & $0.038$ & $0.081$ \\
$256$&$32$&\phantom{0}$8$& $0.028$ & $0.078$ &$0.040$&$0.086$
&$0.051$&$0.106$&$0.053$&$0.111$ & $0.051$ & $0.111$ \\
$256$&$16$&$16$& $0.016$ & $0.063$ &$0.038$&$0.089$
&$0.044$&$0.085$&$0.045$&$0.080$ & $0.033$ & $0.085$ \\
$256$&\phantom{0}$8$&$32$& $0.022$ & $0.068$ &$0.036$&$0.083$
&$0.051$&$0.098$&$0.050$&$0.102$ & $0.051$ & $0.105$ \\
$512$&$64$&\phantom{0}$8$& $0.020$ & $0.073$ &$0.054$&$0.103$
&$0.052$&$0.084$&$0.042$&$0.090$ & $0.039$ & $0.112$ \\
$512$&$32$&$16$& $0.023$ & $0.070$ &$0.046$&$0.083$ &$0.044$ & $0.090$
&$0.049$&$0.092$ & $0.038$ & $0.080$ \\
$512$&$16$&$32$& $0.029$ & $0.067$ &$0.038$&$0.079$ &$0.056$&
$0.098$&$0.052$&$0.099$ & $0.048$ & $0.101$ \\
$512$&\phantom{0}$8$&$64$& $0.025$ & $0.070$ &$0.050$&$0.102$ &$0.047$&
$0.101$&$0.051$&$0.112$ & $0.054$ & $0.105$ \\
\hline
\end{tabular*}
\end{table}
for $\theta\in\{ -0.9, -0.5,0.5,0.9\}$ even for very small samples
sizes. Furthermore, if $T$ gets larger, the results are basically not
affected by the choice of $N$ in these cases. For $\phi=-0.9$, the
nominal level is underestimated for our choice of $T$, but at least if
$T$ grows the approximation of the nominal level becomes more
precise.

%
\begin{table}[b]
\caption{Rejection probabilities of the test (\protect\ref {test})
under the null hypothesis. The data was generated according to model
(\protect\ref{MA})} \label{Table1a}
\begin{tabular*}{\tablewidth}{@{\extracolsep{\fill}}lllllllllll@{}}
\hline
& & & \multicolumn{2}{l}{$\theta=-0.9$} & \multicolumn{2}{l}{$\theta
=-0.5$} & \multicolumn{2}{l}{$\theta=0.5$} &
\multicolumn{2}{l@{}}{$\theta=0.9$} \\[-4pt]
& & & \multicolumn{2}{l}{\hrulefill} & \multicolumn{2}{l}{\hrulefill}
& \multicolumn{2}{l}{\hrulefill} & \multicolumn{2}{l@{}}{\hrulefill} \\
$T$&$N$&$M$& $5 \%$ & $10 \%$ & $5 \%$ & $10 \%$ & $5 \%$ & $10 \%$ &
$5 \%$ & $10 \%$ \\ \hline
\hphantom{0}$64$&\phantom{0}$8$&\phantom{0}$8$ & $0.024$ & $0.073$ & $0.027$ & $0.060$ &$0.045$&$0.091$ &
$0.045$ & $0.096$ \\
$128$&$16$&\phantom{0}$8$ & $0.033$ & $0.071$ & $0.037$ & $0.085$ &$0.043$&$0.087$
& $0.029$ & $0.076$ \\
$128$&\phantom{0}$8$&$16$ & $0.028$ & $0.063$ & $0.031$ & $0.071$ &$0.050$&$0.102$
& $0.028$ & $0.085$ \\
$256$&$32$&\phantom{0}$8$ & $0.047$ & $0.085$ & $0.033$ & $0.081$ &$0.040$&$0.074$
& $0.042$ & $0.080$ \\
$256$&$16$&$16$ & $0.044$ & $0.095$ & $0.031$ & $0.080$
&$0.043$&$0.083$ & $0.035$ & $0.076$ \\
$256$&\phantom{0}$8$&$16$ & $0.029$ & $0.074$ & $0.034$ & $0.081$ &$0.059$&$0.112$
& $0.038$ & $0.076$ \\
$512$&$64$&\phantom{0}$8$ & $0.038$ & $0.084$ & $0.041$ & $0.087$ &$0.052$&$0.106$
& $0.041$ & $0.089$ \\
$512$&$32$&$16$ & $0.047$ & $0.091$ & $0.043$ & $0.073$
&$0.047$&$0.094$ & $0.050$ & $0.100$ \\
$512$&$16$&$32$ & $0.036$ & $0.085$ & $0.044$ & $0.082$
&$0.050$&$0.093$ & $0.050$ & $0.087$ \\
$512$&\phantom{0}$8$&$64$ & $0.051$ & $0.094$ & $0.040$ & $0.078$ &$0.070$&$0.116$
& $0.037$ & $0.080$ \\
\hline
\end{tabular*}
\end{table}

To study the power of the test (\ref{test}), we simulate data from the
following four models
which all correspond to the alternative of non-stationary processes. In
particular, we consider
%
\begin{eqnarray}
\label{alt1} X_{t,T}&=&(1+t/T)Z_t,
\\
\label{alt2} X_{t,T}&=&-0.9\sqrt{\frac{t}{T}}X_{t-1,T}+Z_t,
\\
\label{alt3} X_{t,T}&=& \cases{ %
0.5X_{t-1}+Z_t,
&\quad if $\displaystyle 1\leq t \leq\frac{T}{2}$,
\cr
-0.5X_{t-1}+Z_t,
&\quad if $\displaystyle \frac{T}{2}+1\leq t \leq T$,}
\\
\label{alt4} X_{t,T}&=&Z_t+0.8\cos\bigl(1.5-\cos(4\uppi t/T)
\bigr)Z_{t-q},
\end{eqnarray}
where we display the results for the last model for different $q \in\N
$. Note that due to Remark~\ref{rem3} the alternative (\ref{alt3}) also
fits into the theoretical framework. The corresponding rejection
probabilities are reported in Table \ref{Table2} and we observe a reasonable
behavior of the procedure in the first three considered cases, whereas
power is rather low for the alternative
(\ref{alt4}). Similar to the null
hypothesis we observe robustness with respect to
different choices of $N$, and even for the choice $M=32$, $N=8$,
which appears to be implausible in view of (\ref{assNM}), the results
are satisfying.
It might be of interest to compare these results both with the
pre-periodogram approach from Theorem \ref{thm1a} and with other tests
for the hypothesis of stationarity which have been recently suggested
in the literature. In particular, we consider the tests of Paparoditis
\cite{paparoditis2010}, Dwivedi and Subba~Rao \cite{rao2010} and Dette,
Preuss and Vetter \cite{detprevet2010}.

\begin{table}
\caption{Rejection probabilities of the test (\protect\ref {test}) for
several alternatives} \label{Table2}
\begin{tabular*}{\tablewidth}{@{\extracolsep{\fill}}lllllllllllll@{}}
\hline
&&
& \multicolumn{2}{l}{(\ref{alt1})} & \multicolumn{2}{l}{(\ref{alt2})}
& \multicolumn{2}{l}{(\ref{alt3})} & \multicolumn{2}{l}{(\ref{alt4})
$q=1$} & \multicolumn{2}{l@{}}{(\ref{alt4}) $q=6$} \\[-4pt]
&&
& \multicolumn{2}{l}{\hrulefill} & \multicolumn{2}{l}{\hrulefill}
& \multicolumn{2}{l}{\hrulefill} & \multicolumn{2}{l}{\hrulefill}
& \multicolumn{2}{l@{}}{\hrulefill}\\
$T$&$N$&$M$& $5 \%$ & $10 \%$ & $5 \%$ & $10 \%$ & $5 \%$ & $10 \%$ &
$5 \%$ & $10 \%$ & $5 \%$ & $10 \%$ \\ \hline
\hphantom{0}$64$&\phantom{0}$8$&\phantom{0}$8$&$0.286$&$0.444$&$0.186$&$0.328$&$0.168$&$0.270$& $0.046$ &
$0.098$ & $0.052$ & $0.104$ \\
$128$&$16$&\phantom{0}$8$&$0.686$&$0.772$&$0.396$&$0.546$&$0.308$&$0.466$ &
$0.090$ & $0.154$ & $0.072$ & $0.130$ \\
$128$&\phantom{0}$8$&$16$&$0.624$&$0.758$&$0.382$&$0.578$&$0.410$&$0.548$ &
$0.082$ & $0.144$ & $0.080$ & $0.136$ \\
$256$&$32$&\phantom{0}$8$&$0.958$&$0.974$&$0.672$&$0.814$&$0.742$&$0.912$ &
$0.110$ & $0.186$ & $0.102$ & $0.166$ \\
$256$&$16$&$16$&$0.942$&$0.978$&$0.698$&$0.814$&$0.640$&$0.806$ &
$0.118$ & $0.202$ & $0.098$ & $0.166$ \\
$256$&\phantom{0}$8$&$32$&$0.944$&$0.970$&$0.760$&$0.868$&$0.672$&$0.808$ &
$0.118$ & $0.210$ & $0.086$ & $0.144$ \\
\hline
\end{tabular*}
\end{table}

In Table \ref{Table2b}, we present the rejection frequencies for the
test based on the
pre-periodogram as defined in (\ref{preperio1}). Recall that the use
of the pre-periodogram does not require the specification
of the value $N$, which specifies the number of observations for the
calculation of the local periodogram. This makes its
use attractive for practitioners. However, the results of the
simulation study
show that compared to the local periodogram the use of the
pre-periodogram yields to a substantial loss
of power for
all four alternatives. In particular for alternatives of the form
(\ref{alt3}), the test cannot be recommended.

\begin{table}[b]
\caption{Rejection probabilities of the test based on the
pre-periodogram for several alternatives}
\label{Table2b}
\begin{tabular*}{\tablewidth}{@{\extracolsep{\fill}}lllllllllll@{}}
\hline
& \multicolumn{2}{l}{(\ref{alt1})} & \multicolumn{2}{l}{(\ref{alt2})}
& \multicolumn{2}{l}{(\ref{alt3})} & \multicolumn{2}{l}{(\ref{alt4})
$q=1$} & \multicolumn{2}{l@{}}{(\ref{alt4}) $q=6$} \\[-4pt]
& \multicolumn{2}{l}{\hrulefill} & \multicolumn{2}{l}{\hrulefill}
& \multicolumn{2}{l}{\hrulefill} & \multicolumn{2}{l}{\hrulefill}
& \multicolumn{2}{l@{}}{\hrulefill}\\
$T$& $5 \%$ & $10 \%$ & $5 \%$ & $10 \%$ & $5 \%$
& $10 \%$ & $5 \%$ & $10 \%$ & $5 \%$ & $10 \%$ \\ \hline
\hphantom{0}$64$&$0.188$&$0.340$&$0.080$&$0.202$&$0.022$&$0.056$&$0.024$&$0.076$
&$0.044$&$0.102$ \\
$128$&$0.552$&$0.702$&$0.216$&$0.392$&$0.036$&$0.116$&$0.038$&$0.086$
&$0.052$&$0.098$ \\
$256$&$0.938$&$0.968$&$0.580$&$0.734$&$0.080$&$0.176$&$0.062$&$0.150$
&$0.088$&$0.132$ \\
\hline
\end{tabular*}
\end{table}

In Table \ref{Table2a}, we show the corresponding rejection
probabilities for the test proposed in Dette, Preuss and Vetter
\cite{detprevet2010}, which is the only of the remaining methods
depending on one regularization parameter only. These authors proposed
to estimate the $L_2$ distance
\[
\int_0^1\int_0^\uppi
\biggl(f(u,\lambda)-\int_0^1f(v,\lambda)\,\d v
\biggr)^2 \,\d \lambda \,\d u
\]
using sums of the (squared) periodogram. In order to provide a fair
comparison between the two methods, we also employ the
$\operatorname{AR}(\infty)$
bootstrap to the corresponding test to generate critical values. It
turns out that without a bootstrap the method of Dette, Preuss and
Vetter \cite{detprevet2010} is much more sensitive with respect to
different choices of $N$. We observe that the new method also
outperforms the test proposed by Dette, Preuss and Vetter
\cite{detprevet2010} in the alternatives (\ref{alt1}) and (\ref{alt2}).
In most cases the differences are substantial. On the other hand, for
the alternative (\ref{alt3}) the procedure of Dette, Preuss and Vetter
\cite{detprevet2010} has larger power if $T=64$ and $T=128$, but for
$T=256$ the novel method performs better in this case as well.
Nevertheless, the new approach is clearly outperformed by the proposal
of Dette, Preuss and Vetter \cite{detprevet2010} for the alternative
(\ref{alt4}).

\begin{table}%
\caption{Rejection probabilities of the test proposed by Dette, Preuss
and Vetter \cite{detprevet2010} for several alternatives (quantiles
obtained by $\operatorname{AR}(\infty)$ bootstrap)} \label{Table2a}
\begin{tabular*}{\tablewidth}{@{\extracolsep{\fill}}lllllllllllll@{}}
\hline
&&
& \multicolumn{2}{l}{(\ref{alt1})} & \multicolumn{2}{l}{(\ref{alt2})}
& \multicolumn{2}{l}{(\ref{alt3})} & \multicolumn{2}{l}{(\ref{alt4})
$q=1$} & \multicolumn{2}{l@{}}{(\ref{alt4}) $q=6$} \\[-4pt]
&&
& \multicolumn{2}{l}{\hrulefill} & \multicolumn{2}{l}{\hrulefill}
& \multicolumn{2}{l}{\hrulefill} & \multicolumn{2}{l}{\hrulefill}
& \multicolumn{2}{l@{}}{\hrulefill}\\
$T$&$N$&$M$& $5 \%$ & $10 \%$ & $5 \%$ & $10 \%$ & $5 \%$ & $10 \%$ &
$5 \%$ & $10 \%$ & $5 \%$ & $10 \%$ \\ \hline
\hphantom{0}$64$&\phantom{0}$8$&\phantom{0}$8$&$0.116$&$0.196$&$0.188$&$0.232$&$0.250$&$0.344$& $0.244$ &
$0.350$ & $0.056$ & $0.116$ \\
$128$&$16$&\phantom{0}$8$&$0.106$&$0.160$&$0.256$&$0.330$&$0.370$&$0.552$& $0.490$
& $0.584$ & $0.226$ & $0.336$ \\
$128$&\phantom{0}$8$&$16$&$0.168$&$0.268$&$0.220$&$0.286$&$0.432$&$0.566$& $0.398$
& $0.516$ & $0.072$ & $0.126$ \\
$256$&$32$&\phantom{0}$8$&$0.378$&$0.498$&$0.282$&$0.412$&$0.746$&$0.922$& $0.740$
& $0.836$ & $0.532$ & $0.670$ \\
$256$&$16$&$16$&$0.208$&$0.368$&$0.276$&$0.410$&$0.618$&$0.794$&
$0.716$ & $0.816$ & $0.342$ & $0.444$ \\
$256$&\phantom{0}$8$&$32$&$0.224$&$0.338$&$0.300$&$0.418$&$0.582$&$0.744$& $0.620$
& $0.760$ & $0.104$ & $0.178$ \\ \hline
\end{tabular*}
\end{table}

\begin{table}[b]
\caption{Rejection probabilities of the test proposed by Paparoditis
\cite{paparoditis2010} for several alternatives (quantiles obtained by
$\operatorname{AR}(\infty)$ bootstrap)} \label{Table3a}
\begin{tabular*}{\tablewidth}{@{\extracolsep{\fill}}lllllllllllll@{}}
\hline
&&& \multicolumn{2}{l}{(\ref{alt1})} & \multicolumn{2}{l}{(\ref{alt2})}
& \multicolumn{2}{l}{(\ref{alt3})} & \multicolumn{2}{l}{(\ref{alt4})
$q=1$} & \multicolumn{2}{l@{}}{(\ref{alt4}) $q=6$} \\[-4pt]
&&
& \multicolumn{2}{l}{\hrulefill} & \multicolumn{2}{l}{\hrulefill}
& \multicolumn{2}{l}{\hrulefill} & \multicolumn{2}{l}{\hrulefill}
& \multicolumn{2}{l@{}}{\hrulefill}\\
$T$&$N$&$M$& $5 \%$ & $10 \%$ & $5 \%$ & $10 \%$ & $5 \%$ & $10 \%$ &
$5 \%$ & $10 \%$ & $5 \%$ & $10 \%$ \\ \hline
\hphantom{0}$64$&\phantom{0}$8$&\phantom{0}$8$&$0.054$&$0.126$&$0.050$&$0.122$&$0.078$&$0.170$& $0.058$ &
$0.104$ & $0.034$ & $0.064$ \\
$128$&$16$&\phantom{0}$8$&$0.150$&$0.242$&$0.158$&$0.262$&$0.112$&$0.198$& $0.128$
& $0.218$ & $0.082$ & $0.140$ \\
$128$&\phantom{0}$8$&$16$&$0.066$&$0.154$&$0.120$&$0.254$&$0.166$&$0.270$& $0.080$
& $0.170$ & $0.034$ & $0.066$ \\
$256$&$32$&\phantom{0}$8$&$0.304$&$0.424$&$0.248$&$0.380$&$0.298$&$0.448$& $0.288$
& $0.428$ & $0.102$ & $0.180$ \\
$256$&$16$&$16$&$0.234$&$0.344$&$0.276$&$0.404$&$0.258$&$0.374$&
$0.288$ & $0.420$ & $0.120$ & $0.174$ \\
$256$&\phantom{0}$8$&$32$&$0.126$&$0.226$&$0.240$&$0.374$&$0.298$&$0.376$& $0.158$
& $0.266$ & $0.050$ & $0.106$ \\ \hline
\end{tabular*}
\end{table}

In Table \ref{Table3a}, we show the rejection frequencies for the
method which was proposed in Paparoditis \cite{paparoditis2010}. This concept
basically works by estimating
\[
\sup_{v \in[0,1]} \int_{-\uppi}^\uppi \biggl(
\frac{f(v,\lambda
)}{\int_0^1
f(u,\lambda) \,\d u}-1 \biggr)^2 \,\d \lambda
\]
via a smoothed local periodogram, which requires the choice of a
smoothing bandwidth besides the window length $N$. We choose the
uniform kernel function, and as recommended by the author we select the
bandwidth via the cross validation criterion of Beltr{\~a}o and
Bloomfield \cite{crossvalidation}. To provide a fair comparison, we
also use the $\operatorname{AR}(\infty)$ bootstrap to obtain critical values. For the
alternatives (\ref{alt1})--(\ref{alt3}) the proposal of Paparoditis
\cite{paparoditis2010} yields substantial less power than the approach
proposed in this paper, whereas for the alternative (\ref{alt4}) no
clear picture can be drawn. For $q=1$, the method of Paparoditis
\cite{paparoditis2010} performs better, while there is no significant
difference in the performance if $q=6$. In any case, Paparoditis
\cite{paparoditis2010} is clearly outperformed by the approach of
Dette, Preuss and Vetter \cite{detprevet2010} for (\ref{alt4}).

Finally, we compare our approach to that proposed in Dwivedi and
Subba~Rao \cite{rao2010}. These authors suggested a Portmanteau type
test by estimating
\[
T \sum_{r=1}^m \bigl|c_T(r)\bigr|^2,
\]
where $c_T(r)$ is the covariance of the process at lag $r$. For the
estimation of $c_T(r)$, the authors require the choice of a smoothing
bandwidth, and again we use the cross validation criterion and the
uniform kernel function. Dwivedi and Subba~Rao \cite{rao2010} also have
to choose the maximal lag $m \in N$ up to which they want to estimate
$c_T(r)$, and we pick $m=5$ in the simulations. As in the other
examples, we employ the $\operatorname{AR}(\infty)$ bootstrap, and the results are
given in Table~\ref{Table3b}. A comparison with our method yields a
%
\begin{table}%
\caption{Rejection probabilities of the test proposed by Dwivedi and
Subba~Rao \cite{rao2010} for several alternatives (quantiles obtained
by $\operatorname{AR}(\infty)$ bootstrap)} \label{Table3b}
\begin{tabular*}{\tablewidth}{@{\extracolsep{\fill}}lllllllllll@{}}
\hline
& \multicolumn{2}{l}{(\ref{alt1})} & \multicolumn{2}{l}{(\ref{alt2})}
& \multicolumn{2}{l}{(\ref{alt3})} & \multicolumn{2}{l}{(\ref{alt4})
$q=1$} & \multicolumn{2}{l@{}}{(\ref{alt4}) $q=6$} \\[-4pt]
& \multicolumn{2}{l}{\hrulefill} & \multicolumn{2}{l}{\hrulefill}
& \multicolumn{2}{l}{\hrulefill} & \multicolumn{2}{l}{\hrulefill}
& \multicolumn{2}{l@{}}{\hrulefill}\\
$T$& $5 \%$ & $10 \%$ & $5 \%$ & $10 \%$ & $5 \%$ & $10 \%$ & $5 \%$ &
$10 \%$ & $5 \%$ & $10 \%$ \\ \hline
\hphantom{0}$64$&$0.174$&$0.266$&$0.056$&$0.100$&$0.082$&$0.164$& $0.072$ & $0.120$
& $0.046$ & $0.098$ \\
$128$&$0.274$&$0.386$&$0.058$&$0.114$&$0.122$&$0.208$& $0.126$ &
$0.206$ & $0.092$ & $0.162$ \\
$256$&$0.604$&$0.716$&$0.128$&$0.210$&$0.174$&$0.276$& $0.234$ &
$0.340$ & $0.174$ & $0.272$ \\ \hline
\end{tabular*}
\end{table}
result similar to the approach of Paparoditis \cite{paparoditis2010}.
Our approach performs better for the alternatives
(\ref{alt1})--(\ref{alt3}) while the proposal of Dwivedi and Subba~Rao
\cite{rao2010} yields a higher power in model (\ref{alt4}). Again it is
clearly outperformed in this case by the test proposed in Dette, Preuss
and Vetter \cite{detprevet2010}.

\subsection{Data example}

\begin{figure}

\includegraphics{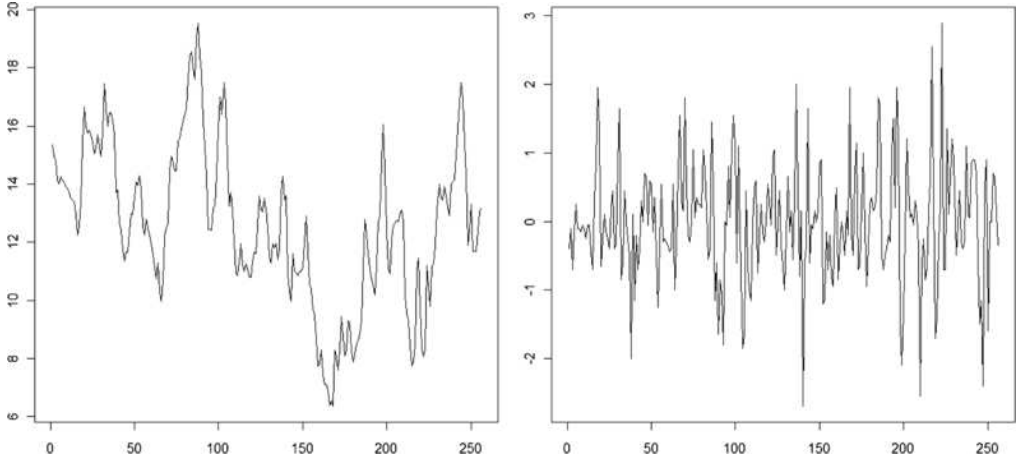}

\caption{Left panel: Weekly egg prices at a
German agriculture market between April 1967 and March~1972. Right
panel: First order difference of the weekly egg prices.}
\label{figure2}
\end{figure}

As an illustration, we consider $T=257$ observations of weekly egg
prices at a German agriculture market between April 1967 and March
1972. A plot of the data is given in Figure \ref{figure2}, and
following Paparoditis \cite{paparoditis2010} the first order difference
$\Delta_t=X_t-X_{t-1}$ of the observed time series are analyzed.
Although several stationary models were proposed in the literature to
fit this data (cf. Paparoditis \cite{paparoditis2010}), the new test
rejects the null hypothesis with $p$-value $0.006$ if we choose $N=32$
or $N=16$, and with $p$-value $0.001$ if we choose $N=8$. These results
demonstrate again that the choice of $N$ does not have too much
influence on the outcome, and that even the somewhat implausible choice
of $N=8$ yields a $p$-value similar to the others.

Note that in Paparoditis \cite{paparoditis2010} a longer version of the above time
series was analyzed, namely $1201$ observations of weekly egg prices
between April 1967 and May 1990. However, we obtain a $p$-value of
exactly $0$ even if we choose $10^6$ bootstrap replicates in this case,
which is why we consider the first $257$ datapoints only. Paparoditis
\cite{paparoditis2010} rejects the null hypothesis of stationarity at level
$5 \%$ if the whole dataset is used, but his approach yields a
$p$-value of $0.1834$ if it is applied to the first $257$ observations
of the time series only, and therefore the hypothesis of stationarity
cannot be rejected at a reasonable size using his method. Roughly the
same $p$-value, namely $0.189$, can be observed if the approach of
Dwivedi and Subba~Rao \cite{rao2010} is employed.

\begin{appendix}\label{app}
\section*{Appendix: Proofs}

\setcounter{subsection}{0}
\subsection{\texorpdfstring{Proof of Theorem \protect\ref{thm1}}{Proof of Theorem 2.1}} \label{sec51}

Throughout the proof, we set $y_j=(v_j,\omega_j) \in[0,1]^2$ for $
j=1,\ldots,K$ and $K \in\N$. To show weak convergence we follow
Theorems 1.5.4 and 1.5.7 in van~der Vaart and Wellner
\cite{wellnervandervaart} and prove the following two claims:
\begin{enumerate}[(2)]
\item[(1)] Convergence of the finite dimensional distributions, that is,
%
\setcounter{equation}{0}
\begin{equation}
\label{finitedistr} \bigl(\hat G_T(y_j)
\bigr)_{j=1,\ldots,K} \xxrightarrow{D} \bigl(G(y_j)
\bigr)_{j=1,\ldots,K}.
\end{equation}
\item[(2)] Stochastic equicontinuity, that is,
%
\begin{equation}
\label{eqcont} \hspace*{-19pt}\forall\eta,\epsilon>0\ \exists\delta>0\dvt\qquad
\lim_{T\rightarrow
\infty} P
\Bigl( \sup_{y_1, y_2 \in[0,1]^2\dvtx d_2(y_1,y_2) <\delta} \bigl|G_T(y_1)-G_T(y_2)\bigr|
> \eta \Bigr)< \epsilon,
\end{equation}
where $d_2(y_1,y_2)=\sqrt{(v_1-v_2)^2+(w_1-w_2)^2}$.
\end{enumerate}
\begin{pf*}{Proof of (\ref{finitedistr})}
The claim follows from similar arguments as given in the proof of
Theorem~3.1 in Dette, Preuss and Vetter \cite{detprevet2010}. For the
sake of brevity and because we will use similar arguments in the proof
of (\ref{eqcont}), we will sketch how the assertions
%
\begin{eqnarray}
\label{expect}
\hspace*{-18pt}\E\bigl(\hat G_T(v,\omega)\bigr)
&\xrightarrow{T
\rightarrow\infty}& 0,
\\
\label{cov}
\hspace*{-18pt}\Cov\bigl(\hat G_T(y_1), \hat
G_T(y_2)\bigr)
&\xrightarrow{T \rightarrow \infty}&
\frac{1}{2\uppi}\int_0^1\int
_0^{\uppi\min(\omega_1,\omega
_2)} \bigl(1_{[0,v_1]}(u)-v_1
\bigr) \bigl(1_{[0,v_2]}(u)-v_2\bigr)
\nonumber\\[-8pt]\\[-8pt]
\hspace*{-18pt}&&\hspace*{83.4pt}
{}\times f^2(u,\lambda) \,\d \lambda \,\d u\nonumber
\end{eqnarray}
can be shown. Note that we have
\begin{eqnarray*}
\hat G_T(v,\omega)&=&\frac{1}{\sqrt{T}}\sum
_{j=1}^M\sum_{k=1}^{
{N}/{2}}
\phi_{v,\omega,M,N}(u_j,\lambda_k)I_N^X(u_j,
\lambda_k)-\sqrt{T} D_{N,M}(\phi_{v,\omega,M,N})
\\
&=&\!:G_T(\phi_{v,\omega,M,N})
\end{eqnarray*}
with
\[
\phi_{v,\omega,M,N}(u,\lambda):=\biggl(I_{[0,{\lfloor vM \rfloor
}/{M}]}(u)-\frac{\lfloor vM \rfloor}{M}
\biggr)I_{[0,{ 2 \uppi\lfloor
\omega{N}/{2} \rfloor/N}]}(\lambda)
\]
for $u,\lambda\geq0$ and
\[
D_{N,M}(\phi):=\frac{1}{2\uppi}\int_0^1
\int_0^\uppi\phi(u,\lambda) f(u,\lambda) \,\d \lambda
\,\d u.
\]
In order to simplify some technical arguments, we also define
\[
\phi_{v,\omega,M,N}(u,\lambda):=\phi_{v,\omega,M,N}(u,-\lambda)
\]
for $ u\geq0, \lambda<0$
and obtain from (\ref{apprbed})
\begin{eqnarray*}
&& \E \Biggl( \frac{1}{T}\sum_{j=1}^M
\sum_{k=1}^{{N}/{2}} \phi_{v,\omega,M,N}(u_j,
\lambda_k)I_N^X(u_j,
\lambda_k) \Biggr)
\\
&&\quad= \frac{1}{T}\sum_{j=1}^M\sum
_{k=1}^{{N}/{2}} \phi_{v,\omega,M,N}(u_j,
\lambda_k) \\
&&\hspace*{39.6pt}\qquad{}\times\frac{1}{2\uppi N}\sum_{p,q=0}^{N-1}
\sum_{l,m=-\infty}^\infty\psi_l \biggl(
\frac{t_j-N/2+1+p}{T} \biggr)\psi_m \biggl(\frac{t_j-N/2+1+q}{T} \biggr)
\\
&&\hspace*{151.3pt}{}\times \E (Z_{t_j-N/2+1+p-m}Z_{t_j-N/2+1+q-l})\\
&&\hspace*{151.3pt}{}\times \exp\bigl(-\i \lambda_k(p-q)
\bigr) \bigl( 1 +\O(1/T) \bigr).
\end{eqnarray*}
A Taylor expansion now yields that this term becomes
\begin{eqnarray*}
&& \frac{1}{T}\sum_{j=1}^M\sum
_{k=1}^{{N}/{2}} \phi_{v,\omega,M,N}(u_j,
\lambda_k) \frac{1}{2\uppi N}\\[-2pt]
&&\hspace*{39pt}{}\times\sum_{p,q=0}^{N-1}
\sum_{l,m=-\infty}^\infty\psi_l(u_j)
\psi_m(u_j)\\[-2pt]
&&\hspace*{103pt}{}\times
\E(Z_{t_j-N/2+1+p-m}Z_{t_j-N/2+1+q-l})\\[-2pt]
&&\hspace*{103pt}{}\times \exp\bigl(-\i \lambda_k(p-q)
\bigr) \bigl(1+\O(1/T)+\O\bigl(N^2/T^2\bigr)\bigr).
\end{eqnarray*}
See Dette, Preuss and Vetter \cite{detprevet2010} for details. Since
$\E(Z_iZ_j)=0$ for $i \not= j$, we obtain the equation $p=q+m-l$ which
shows that the above expression equals
\begin{eqnarray*}
&& \frac{1}{2\uppi NT}\sum_{j=1}^M\sum
_{k=1}^{{N}/{2}} \phi_{v,\omega,M,N}(u_j,
\lambda_k) \\[-2pt]
&&\qquad\hspace*{39pt}{}\times\sum_{l,m=-\infty}^\infty
\mathop{\sum_{q=0}}_{0 \leq q+m-l \leq N-1}^{N-1}
\psi_l(u_j)\psi_m(u_j) \exp
\bigl(-\i \lambda_k(m-l)\bigr)
\\[-2pt]
&&\qquad{}+\O(1/T)+\O\bigl(N^2/T^2\bigr)
\\[-2pt]
&&\quad= \frac{1}{2\uppi NT}\sum_{j=1}^M\sum
_{k=1}^{{N}/{2}} \phi_{v,\omega,M,N}(u_j,
\lambda_k) \\[-2pt]
&&\qquad\hspace*{59.3pt}{}\times\mathop{\sum_{l,m=-\infty}}_{|l-m| \leq N-1}^\infty
\mathop{\sum_{q=0}}_{0 \leq q+m-l \leq N-1}^{N-1}
\psi_l(u_j)\psi_m(u_j) \exp
\bigl(-\i \lambda_k(m-l)\bigr)
\\[-2pt]
&&\qquad{} + \frac{1}{2\uppi NT}\sum_{j=1}^M\sum
_{k=1}^{{N}/{2}} \phi_{v,\omega,M,N}(u_j,
\lambda_k) \\[-2pt]
&&\qquad\hspace*{71.5pt}{}\times\mathop{\sum_{l,m=-\infty}}_{|l-m| \geq N}^\infty
\mathop{\sum_{q=0}}_{0 \leq q+m-l \leq N-1}^{N-1}
\psi_l(u_j)\psi_m(u_j) \exp
\bigl(-\i \lambda_k(m-l)\bigr)
\\[-2pt]
&&\qquad{} + \O(1/T)+\O\bigl(N^2/T^2\bigr).
\end{eqnarray*}
Dropping the extra condition $0 \leq q+m-l \leq N-1$, the second term
is bounded by
%
\begin{eqnarray}
\label{boundgrN}
C \mathop{\sum_{l,m=-\infty}}_{|l-m| \geq N}^\infty
\sup_u \bigl|\psi_l(u)\bigr| \sup_u\bigl|
\psi_m(u)\bigr|
&\leq& 2C \sum_{m=-\infty}^\infty
\sup_u\bigl|\psi_m(u)\bigr| \mathop{\sum
_{l=-\infty}}_{|l| \geq N/2}^\infty\sup_u \bigl|
\psi_l(u)\bigr|
\nonumber\\[-2pt]
&\leq& \frac{4C \sum_{m=-\infty}^\infty\sup_u|\psi_m(u)|
\sum_{l=-\infty}^\infty|l|\sup_u|\psi_l(u)|}{N}\hspace*{29pt}
\\[-2pt]
&=& \O(1/N)\nonumber
\end{eqnarray}
for some $C \in\R$ and the order follows from (\ref{2.1a}). Using
(\ref{2.1a}) and (\ref{boundgrN}) in the same way again, the first
quantity above can be shown to be equal to
\[
\frac{1}{2\uppi T}\sum_{j=1}^M\sum
_{k=1}^{{N}/{2}} \phi_{v,\omega,M,N}(u_j,
\lambda_k) \sum_{l,m=-\infty}^\infty
\psi_l(u_j)\psi_m(u_j) \exp
\bigl(-\i \lambda_k(m-l)\bigr)+\O(1/N),
\]
and therefore we obtain
\begin{eqnarray*}
&& \E \Biggl( \frac{1}{T}\sum_{j=1}^M
\sum_{k=1}^{{N}/{2}} \phi_{v,\omega,M,N}(u_j,
\lambda_k)I_N^X(u_j,
\lambda_k) \Biggr)
\\
&&\quad= \frac{1}{T}\sum_{j=1}^M\sum
_{k=1}^{{N}/{2}} \phi_{v,\omega,M,N}(u_j,
\lambda_k) f(u_j,\lambda_k)+\O(1/N)\\
&&\qquad{}
+\O\bigl(N^2/T^2\bigr)+\O(1/T)
\\
&&\quad= D_{N,M}(\phi_{v,\omega,M,N}) +\O(1/N)+\O\bigl(N^2/T^2
\bigr)+\O(1/T),
\end{eqnarray*}
where the order of the Riemann approximation follows from the specific
choice of the midpoints~$u_j$. This together with (\ref{assNM}) yields
(\ref{expect}).

To prove (\ref{cov}), we use symmetry arguments and obtain
\begin{eqnarray*}
\hspace*{-4pt}&&T \cum\Biggl(\frac{1}{T}\sum_{j_1=1}^M
\sum_{k_1=1}^{{N}/{2}} \phi_{v_1,\omega_1,M,N}(u_{j_1},
\lambda_{k_1})I_N^X(u_{j_1},
\lambda_{k_1}),
\\
\hspace*{-4pt}&&\hspace*{33.3pt}\frac{1}{T}\sum_{j_2=1}^M
\sum_{k_2=1}^{{N}/{2}} \phi_{v_2,\omega_2,M,N}(u_{j_2},
\lambda_{k_2})I_N^X(u_{j_2},
\lambda_{k_2})\Biggr)
\\
\hspace*{-4pt}&&\quad=\frac{1}{4T}\frac{1}{(2\uppi N)^2}
\sum_{j_1,j_2=1}^M
\sum_{k_1,k_2=-\lfloor({N-1})/{2} \rfloor}^{{N}/{2}} \phi_{v_1,\omega
_1,M,N}(u_{j_1},
\lambda_{k_1})\phi_{v_2,\omega
_2,M,N}(u_{j_2},
\lambda_{k_2})
\\
\hspace*{-4pt}&&\qquad{} \times \sum_{p_1,p_2,q_1,q_2=0}^{N-1} \sum
_{m_1,m_2,l_1,l_2=-\infty}^\infty \psi_{m_1}(u_{j_1})
\psi_{l_1}(u_{j_1})\psi_{m_2}(u_{j_2})
\psi_{l_2}(u_{j_2})
\\
\hspace*{-4pt}&&\hspace*{116.5pt}\qquad{} \times\exp\bigl(-\i \lambda_{k_1}(p_1-q_1)\bigr)
\exp\bigl(-\i \lambda_{k_2}(p_2-q_2)\bigr)
\\
\hspace*{-4pt}&&\hspace*{116.5pt}\qquad{} \times
\cum (Z_{t_{j_1}-N/2+1+p_1-m_1}Z_{t_{j_1}-N/2+1+q_1-l_1},\\
\hspace*{-4pt}&&\hspace*{149.7pt}\qquad{} Z_{t_{j_2}-N/2+1+p_2-m_2}Z_{t_{j_2}-N/2+1+q_2-l_2})
\\
\hspace*{-4pt}&&\hspace*{116.5pt}\qquad{} \times\bigl(1 +\O\bigl(N^2/T^2\bigr)+\O(1/T)\bigr)
\end{eqnarray*}
in the same way as above. Because of
\begin{eqnarray*}
&&\hspace*{-5pt}\cum (Z_{t_{j_1}-N/2+1+p_1-m_1}Z_{t_{j_1}-N/2+1+q_1-l_1},Z_{t_{j_2}-N/2+1+p_2-m_2}Z_{t_{j_2}-N/2+1+q_2-l_2})
\\
&&\hspace*{-5pt}\!\!\quad=\cum(Z_{t_{j_1}-N/2+1+p_1-m_1}Z_{t_{j_2}-N/2+1+q_2-l_2})\cum (Z_{t_{j_2}-N/2+1+p_2-m_2}Z_{t_{j_1}-N/2+1+q_1-l_1})
\\
&&\hspace*{-5pt}\!\!\qquad{}+\cum(Z_{t_{j_1}-N/2+1+p_1-m_1}Z_{t_{j_2}-N/2+1+p_2-m_2})\cum (Z_{t_{j_1}-N/2+1+q_1-l_1}Z_{t_{j_2}-N/2+1+q_2-l_2}),
\end{eqnarray*}
the calculation of the highest order term in the variance splits into
two sums and we only consider the first one (the second sum is treated
completely analogously), which equals
\begin{eqnarray*}
&&\frac{1}{4T}\sum_{j_1,j_2=1}^M\sum
_{k_1,k_2=-\lfloor({N-1})/{2}
\rfloor}^{{N}/{2}} \phi_{v_1,\omega_1,M,N}(u_{j_1},
\lambda_{k_1})\phi_{v_2,\omega_2,M,N}(u_{j_2},
\lambda_{k_2})\frac{1}{(2\uppi N)^2}
\\
&&\qquad{} \times  \sum_{m_1,m_2,l_1,l_2=-\infty}^\infty
\mathop{\mathop{\sum_{q_1,q_2=0}}_{0 \leq
q_2+m_1-l_2+t_{j_2}-t_{j_1} \leq
N-1}}_{
0 \leq q_1+m_2-l_1+t_{j_1}-t_{j_2} \leq N-1}^{N-1}
\psi_{m_1}(u_{j_1})\psi_{l_1}(u_{j_1})
\psi_{m_2}(u_{j_2})\psi_{l_2}(u_{j_2})
\\
&&\hspace*{184pt}{} \times \exp\bigl(-\i (\lambda_{k_1}-\lambda_{k_2})
(q_2-q_1+t_{j_2}-t_{j_1})\bigr)\\
&&\hspace*{184pt}{} \times \exp
\bigl(-\i \lambda_{k_1}(m_1-l_2)-\i
\lambda_{k_2}(m_2-l_1)\bigr)
\\
&&\quad=\frac{1}{4T}\sum_{j_1,j_2=1}^M\sum
_{k_1,k_2=-\lfloor({N-1})/{2}
\rfloor}^{{N}/{2}} \phi_{v_1,\omega_1,M,N}(u_{j_1},
\lambda_{k_1})\phi_{v_2,\omega_2,M,N}(u_{j_2},
\lambda_{k_2})\frac{1}{(2\uppi N)^2}
\\
&&\qquad{} \times
\mathop{\sum_{m_1,m_2,l_1,l_2=-\infty}}_{(+)}^\infty
\mathop{\mathop{\sum_{q_1,q_2=0}}_{0 \leq
q_2+m_1-l_2+t_{j_2}-t_{j_1} \leq N-1}}_{0 \leq
q_1+m_2-l_1+t_{j_1}-t_{j_2} \leq N-1}^{N-1}
\psi_{m_1}(u_{j_1})\psi_{l_1}(u_{j_1})
\psi_{m_2}(u_{j_2})\psi_{l_2}(u_{j_2})
\\
&&\hspace*{184pt}{} \times\exp\bigl(-\i (\lambda_{k_1}-\lambda_{k_2})
(q_2-q_1+t_{j_2}-t_{j_1})\bigr)\\
&&\hspace*{184pt}{} \times
\exp\bigl(-\i \lambda_{k_1}(m_1-l_2)-\i
\lambda_{k_2}(m_2-l_1)\bigr)
\\
&&\hspace*{184pt}{} \times\bigl(1 +\O(1/N)\bigr),
\end{eqnarray*}
where $\sum_{(+)}$ means that summation is only performed over
those indices
$x,y \in\{m_1,m_2,l_1,l_2 \}$ such that $|x-y|<N$,
and the $\O(1/N)$-term follows with (\ref{boundgrN}). Assume that $j_1$
has been chosen. Then $j_2$ must be equal to $j_1,j_1-1$ or $j_1+1$, as
all other combination of $j_1$ and $j_2$ vanish, because of the
condition $0 \leq q_2+m_1-l_2+t_{j_2}-t_{j_1} \leq N-1$
and the fact that the summation is only performed with respect to the
indices satisfying $|x-y|<N$.
If $j_2$ equals $j_1-1$ or $j_1+1$, it follows from the conditions on
$q_1$ and $q_2$ that for chosen $m_i$ and $l_i$, there are at most
$|m_2-l_1|$ possible choices for $q_1$ and at most $|m_1-l_2|$ possible
choices for $q_2$. It therefore follows with (\ref{2.1a}) that the
terms with $j_2 \in\{j_1-1,j_1+1\}$ are of order $\O(1/N)$.

Therefore, we only have to consider the case $j_1=j_2$, and the above
expression is
%
\begin{eqnarray}
\label{varpr1} &&\frac{1}{4T}\sum_{j_1=1}^M
\sum_{k_1,k_2=-\lfloor({N-1})/{2}
\rfloor
}^{{N}/{2}} \phi_{v_1,\omega_1,M,N}(u_{j_1},
\lambda_{k_1})\phi_{v_2,\omega_2,M,N}(u_{j_1},
\lambda_{k_2})
\nonumber\\
&&\quad{} \times \frac{1}{(2\uppi N)^2} \mathop{\sum_{m_1,m_2,l_1,l_2=-\infty}}_{(+)}^\infty
\mathop{\mathop{\sum_{q_1,q_2=0}}_{0 \leq q_2+m_1-l_2 \leq N-1}}_{0
\leq q_1+m_2-l_1 \leq N-1}^{N-1}
\psi_{m_1}(u_{j_1})\psi_{l_1}(u_{j_1})
\psi_{m_2}(u_{j_2})\psi_{l_2}(u_{j_2})
\nonumber\\
&&\hspace*{171pt}\quad{} \times
\exp\bigl(-\i (\lambda_{k_1}-\lambda_{k_2})
(q_2-q_1)\bigr)\\
&&\hspace*{171pt}\quad{}\times\exp\bigl(-\i \lambda_{k_1}(m_1-l_2)-\i
\lambda_{k_2}(m_2-l_1)\bigr) \hspace*{27pt}\nonumber\\
&&\hspace*{171pt}\quad{}\times\bigl(1+\O(1/N)
\bigr).\nonumber
\end{eqnarray}
Observing
\[
\frac{1}{N} \sum_{q=0}^{N-1} \exp
\bigl(-\i  (\lambda_{k_1}-\lambda_{k_2}) q\bigr) = %
\cases{ 1, &\quad $k_1-k_2=lN$ with $l \in\mathbb{Z}$,
\cr
0, &\quad
else,}
\]
it follows that for fixed $m_1$, $l_2$ and $k_1 \not= k_2$ we have
\begin{eqnarray*}
\Biggl| \mathop{\sum_{q_2=0}}_{0 \leq q_2+m_1-l_2 \leq N-1}^{N-1}
\exp\bigl(-\i  (\lambda_{k_1}-\lambda_{k_2})q_2\bigr)
\Biggr| &=& \Biggl| \mathop{\mathop{\mathop{\sum_{q_2=0}}_{q_2+m_1-l_2 <0}}_{\mathrm{or}}}_{q_2+m_1-l_2 > N-1}^{N-1}
\exp\bigl(-\i  (\lambda_{k_1}-\lambda_{k_2})q_2\bigr)
\Biggr|
\\
&\leq& |m_1-l_2|,
\end{eqnarray*}
which implies
%
\begin{equation}
\label{k1unglk2}
\hspace*{-12pt}\Biggl| \frac{1}{(2\uppi N)^2} \mathop{\mathop{\sum
_{q_1,q_2=0}}_{0 \leq q_2+m_1-l_2 \leq N-1}}_{0
\leq q_1+m_2-l_1 \leq N-1}^{N-1} \exp
\bigl(-\i (\lambda_{k_1}-\lambda_{k_2}) (q_2-q_1)
\bigr) \Biggr| \leq |m_1-l_2||m_2-l_1|/(2
\uppi N)^2.
\end{equation}
By using (\ref{2.1a}) and (\ref{k1unglk2}), it can now be seen that
all terms
with $k_1 \not= k_2$ are of the order $\O(1/N)$, and
similar arguments as used in the calculation of the expectation yield
that (\ref{varpr1}) equals
\begin{eqnarray*}
&&
\frac{1}{4\uppi}\int_0^1\int
_0^{\uppi\min(\omega_1,\omega_2)} \bigl(1_{[0,v_1]}(u)-v_1
\bigr) \bigl(1_{[0,v_2]}(u)-v_2\bigr)f^2(u,\lambda) \,\d
\lambda \,\d u\\
&&\quad{}+\O(1/N)+\O\bigl(N^2/T^2\bigr).
\end{eqnarray*}
\upqed\end{pf*}
\begin{pf*}{Proof of (\ref{eqcont})}
Note that
\[
\mathcal{F}_T:= \bigl\{\phi_{v,\omega,M,N};v,\omega\in[0,1] \bigr\}
= \bigl\{\phi_{v,\omega,M,N}; (v,\omega) \in P_T \bigr\},
\]
where
\[
P_T:= \biggl\{0,\frac{1}{M},\frac{2}{M},\ldots,
\frac{M-1}{M},1 \biggr\} \times \biggl\{0,\frac{2}{N},
\frac{4}{N},\ldots,1-\frac{2}{N},1 \biggr\} %
\]
(recall that $N$ is assumed to be even throughout this paper).
We define
\[
\rho_2(\phi):= \biggl( \int_0^1
\int_0^\uppi\phi^2(u,\lambda) \,\d \lambda
\,\d u \biggr)^{1/2},
\]
and $\mathcal{F}^2_T$ is the set of functions, which can be expressed
as a sum or a difference of two elements in $\mathcal{F}_T$. The main
task is to prove the following theorem.
%
\begin{satz}
\label{appth1}
There exists a constant $C \in\R$ such that for all $\phi\in
\mathcal{F}^2_T$:
\[
\E\bigl(\bigl|\hat G_T(\phi)\bigr|^k\bigr) \leq(2k)!
C^k \rho_2(\phi)^k\qquad \forall k \in\N\mbox{ even}.
\]
\end{satz}
Stochastic equicontinuity follows then by similar arguments as given in
Dahlhaus \cite{dahlhaus1988}, which is why we will only sketch the main steps
and refer to his work for most details.
The first consequence of Theorem \ref{appth1} regards the existence of
a constant $C_1 \in\R$ such that for all $g,h \in\mathcal{F}_T$ and
$\eta>0$:
\[
P\bigl(\bigl|\hat G_T(g)-\hat G_T(h)\bigr|>\eta
\rho_2(g-h)\bigr) \leq96 \exp\biggl(-\sqrt {\frac
{\eta}{C_1}}
\biggr).
\]
A straightforward modification of the chaining lemma in Chapter VII.2
of Pollard \cite{pollard} then yields that for a stochastic process $(Z(v))_{v
\in V}$, whose index set $V$ has a finite covering integral
%
\begin{equation}
\label{covering} J(\delta)=\int_0^\delta \biggl[
\log \biggl(\frac{48N(u)^2}{u} \biggr) \biggr]^2 \,\d u
\end{equation}
for all $\delta$ and which satisfies
\[
P \bigl( \bigl|Z(v)-Z(w)\bigr|>\nu d(v,w) \bigr) \leq96 \exp\biggl(-\sqrt{
\frac
{\nu}{C_1}}\biggr)
\]
for a semi-metric $d$ on $V$ and a constant $C_1 \in\R$, there exist a
dense subset $V^* \subset V$
such that
\[
P \bigl( \exists v,w \in V^* \mbox{ with } d(v,w)<\epsilon\mbox{ and }
\bigl|Z(v)-Z(w)\bigr|>26C_1J\bigl(d(v,w)\bigr) \bigr) \leq2 \epsilon.
\]
In (\ref{covering}), $N(u)$ is the covering number which is defined as
the smallest number $m \in\N$ for which
there exist $z_1,\ldots,z_m \in V$ with $
\min_{i}d(z,z_i) \leq u $ for all $z \in V$.
By using $y_i=(v_i, \omega_i)$, we obtain
\begin{eqnarray*}
&&P \Bigl( \sup_{y_1,y_2 \in P_T\dvtx d_2(y_1,y_2) < \delta} \bigl|\hat G_T(v_2,w_2)-
\hat G_T(v_1,w_1)\bigr|>\eta \Bigr)
\\
&&\quad\leq P \Bigl( \sup_{f,g \in\mathcal{F}_T\dvtx
\rho_2(f,g) < \epsilon
(\delta)} \bigl|\hat G_T(f)-\hat
G_T(g)\bigr|>\eta \Bigr)
\end{eqnarray*}
for a certain sequence $\epsilon(\delta) \xrightarrow{\delta
\rightarrow0} 0$ by continuity.
The right-hand side of this inequality equals
\begin{eqnarray*}
&& P \Bigl( \sup_{f,g \in\mathcal{F}_T\dvtx \rho_2(f,g) < \epsilon
(\delta)} \bigl|\hat G_T(f)-\hat
G_T(g)\bigr|>\eta, \eta\geq26C_1J_T\bigl(\epsilon(
\delta)\bigr) \Bigr)
\\
&&\qquad{}+P \Bigl( \sup_{f,g \in\mathcal{F}_T\dvtx \rho_2(f,g) < \epsilon
(\delta)} \bigl|\hat G_T(f)-\hat
G_T(g)\bigr|>\eta, \eta< 26C_1J_T\bigl(\epsilon(
\delta)\bigr) \Bigr)
\\
&&\quad\leq 2\epsilon(\delta)+P\bigl(\eta< 26C_1J_T\bigl(
\epsilon(\delta)\bigr)\bigr),
\end{eqnarray*}
where $J_T(\delta)$ is the corresponding covering integral of
$\mathcal
{F}_T$. Note that $\eta< 26C_1J_T(\epsilon(\delta))$ is not random and
that $J_T(\delta)$ can be bounded by $J(\delta)$, which is the covering
integral of $\bigcup_{i=1}^{\infty} \mathcal{F}_i$ (which is finite for
every $\delta$). Because of $J(\epsilon(\delta))\xrightarrow{\delta
\rightarrow0}0$, we have
$\eta> 26C_1J(\delta)$ whenever $\delta$ is sufficiently small and obtain
\[
P \Bigl( \sup_{f,g \in\mathcal{F}_T\dvtx \rho_2(f,g)
< \epsilon(\delta)} \bigl|\hat G_T(f)-\hat G_T(g)\bigr|>
\eta \Bigr) < 2\epsilon(\delta),
\]
which implies the stochastic equicontinuity.
\begin{pf*}{Proof of Theorem \ref{appth1}} We show
%
\begin{equation}
\label{appineq} \bigl|\cum_l\bigl(\sqrt{T} \hat D_T(\phi)
\bigr)\bigr| \leq(2l)! \tilde C^l \rho_2(\phi)^l\qquad
\forall l \in\N,
\end{equation}
where
\[
\hat D_T(\phi):= \frac{1}{\sqrt{T}}\hat G_T(\phi)+
D_{N,M}(\phi).
\]
Since $D_{N,M}(\phi)$ is constant, this implies
\[
\bigl|\cum_l(\hat G_T)\bigr| \leq(2l)! C^l
\rho_2(\phi)^l\qquad\forall l \in \N
\]
for some $C$, and then it follows as in Dahlhaus \cite{dahlhaus1988} that
\begin{eqnarray*}
\E\bigl(\bigl|\hat G_T(\phi)\bigr|^k\bigr)&=& \Biggl| \mathop{\mathop{
\sum_{\{P_1,\ldots,P_m\}}}_{\mathrm{Partition}\ \mathrm{of}}}_{\{
1,\ldots,k\}} \Biggl\{
\prod_{j=1}^m \cum_{|P_j|}\bigl(
\hat G_T(\phi)\bigr) \Biggr\} \Biggr| \\
&\leq&\rho_2(
\phi)^k C^k \mathop{\mathop{\sum
_{\{P_1,\ldots,P_m\}}}_{\mathrm{Partition}\ \mathrm{of}}}_{\{
1,\ldots,k\}} \prod
_{j=1}^m \bigl(2|P_j|\bigr)!
\\
&\leq&(2k)! C^k 2^k \rho_2(
\phi)^k,
\end{eqnarray*}
since we only consider the case where $k$ is even. This yields the
assertion.

In order\vspace*{1pt} to prove (\ref{appineq}), we assume without loss of generality
that $l$ is even, as the case for odd $l$ is proved in the same way.
The $l$th cumulant of $\sqrt{T} \hat D_T(\phi)$ is given by
\begin{eqnarray*}
\hspace*{-4pt}&&\frac{1}{2^lT^{l/2}}\sum_{j_1,\ldots,j_l=1}^M\sum
_{k_1,\ldots,k_l=-\lfloor
({N-1})/{2} \rfloor}^{{N}/{2}} \phi(u_{j_1},
\lambda_{k_1})\cdots \phi(u_{j_l},\lambda_{k_l})\frac{1}{(2\uppi N)^l}
\\
\hspace*{-4pt}&&\quad{} \times\sum_{p_1,q_1,p_2,\ldots,p_l,q_l=0}^{N-1}
\sum_{m_1,n_1,m_2,\ldots,m_l,n_l=-\infty}^\infty \psi_{m_1}(u_{j_1})
\cdots\psi_{n_l}(u_{j_l})
\\
\hspace*{-4pt}&&\hspace*{158pt}\quad{} \times\cum
(Z_{t_{j_1}-N/2+1+p_1-m_1}Z_{t_{j_1}-N/2+1+q_1-n_1},\ldots,\\
\hspace*{-4pt}&&\hspace*{158pt}\hspace*{33pt}\quad{} Z_{t_{j_l}-N/2+1+p_l-m_l}Z_{t_{j_l}-N/2+1+q_l-n_l})
\\
\hspace*{-4pt}&&\hspace*{158pt}\quad{} \times\exp\bigl(-\i \lambda_{k_1}(p_1-q_1)
\bigr)\cdots\exp\bigl(-\i \lambda_{k_l}(p_l-q_l)
\bigr) \\
\hspace*{-4pt}&&\hspace*{158pt}\quad{} \times
\bigl(1+\O\bigl(N^2/T^2\bigr)+\O(1/T)\bigr),
\end{eqnarray*}
where both $\O(\cdot)$-terms follow as in the proof of (\ref{expect}).
We define $Y_{i,1}:=Z_{t_{j_i}-N/2+1+p_i-m_i}$ and
$Y_{i,2}:=Z_{t_{j_i}-N/2+1+q_i-n_i}$ for $i \in\{1,\ldots,l\}$. Theorem
2.3.2 in Brillinger \cite{brillinger1981} yields
\[
\cum_l\bigl(\sqrt{T} \hat D_T(\phi)\bigr)=\sum
_\nu V_T(\nu) \bigl(1+\O\bigl(N^2/T^2
\bigr)+\O(1/T)\bigr),
\]
where the sum runs over all indecomposable partitions $\nu= \nu_1
\cup\cdots\cup\nu_{l}$ with $|\nu_i|=2$ ($1 \leq i \leq l$, due to
Gaussianity) of the matrix
%
\begin{equation}
\label{part1} %
\matrix{ Y_{1,1} & Y_{1,2}
\cr
\vdots& \vdots
\cr
Y_{l,1} & Y_{l,2}} %
\end{equation}
and
\begin{eqnarray*}
V_T(\nu)&:=&\frac{1}{2^lT^{l/2}}\sum_{j_1,\ldots,j_l=1}^M
\sum_{k_1,\ldots,k_l=-\lfloor({N-1})/{2} \rfloor}^{{N}/{2}} \phi (u_{j_1},
\lambda_{k_1})\cdots\phi(u_{j_l},\lambda_{k_l})
\\
&&{} \times\frac{1}{(2\uppi N)^l}\sum_{p_1,\ldots,q_l=0}^{N-1}
\sum_{m_1,\ldots,n_l=-\infty}^\infty\psi_{m_1}(u_{j_1})
\cdots\psi_{n_l}(u_{j_l})
\\
&&\hspace*{129pt}{} \times \cum\bigl(Y_{i,k};(i,k) \in\nu_1\bigr) \cdots\cum
\bigl(Y_{i,k};(i,k) \in\nu_l\bigr)
\\
&&\hspace*{129pt}{} \times\exp\bigl(-\i \lambda_{k_1}(p_1-q_1)\bigr)
\cdots\exp\bigl(-\i \lambda_{k_l}(p_l-q_l)
\bigr).
\end{eqnarray*}
We now fix one indecomposable partition $\tilde\nu$ and assume without
loss of generality that
\[
\tilde\nu=\bigcup_{i=1}^{l-1}(Y_{i,1},Y_{i+1,2})
\cup(Y_{l,1},Y_{1,2}). %
\]
Because of $\cum(Z_i,Z_j)\not=0$ for $i \not= j$, we obtain the
following $l$ equations:
%
\begin{eqnarray}
\label{einschr0a}
q_1&=&p_l+n_1-m_l+t_{j_l}-t_{j_1},
\\
\label{einschr0b} q_{i+1}&=&p_{i}+n_{i+1}-m_{i}+t_{j_{i}}-t_{j_{i+1}}
\qquad\mbox{for } i \in\{1,\ldots,l-1\}
\end{eqnarray}
and therefore only $l$ variables (namely $p_i$ for $i \in\{1,\ldots,l\}$)
of the $2l$ variables $p_1,q_1,p_2,\ldots,q_l$ are free to choose and must
satisfy the following conditions:
%
\begin{eqnarray}
\label{einschr1a} 0 &\leq& p_i +n_{i+1}-m_i+t_{j_i}-t_{j_{i+1}}
\leq N-1 \qquad\mbox{for } i \in\{1,\ldots,l-1\},
\\
\label{einschr1b} 0 &\leq&p_l+n_1-m_l+t_{j_l}-t_{j_1}
\leq N-1.
\end{eqnarray}
Using the identities (\ref{einschr0a}) and (\ref{einschr0b}), we obtain
that $V_T(\tilde\nu)$ equals
\begin{eqnarray*}
&& \frac{1}{2^lT^{l/2}}\sum_{j_1,\ldots,j_l=1}^M\sum
_{k_1,\ldots,k_l=-\lfloor
({N-1})/{2} \rfloor}^{{N}/{2}} \phi(u_{j_1},
\lambda_{k_1})\cdots \phi(u_{j_l},\lambda_{k_l})
\frac{1}{(2\uppi N)^l}\\
&& \hspace*{0pt}\quad{}\times\sum_{p_1,p_2,\ldots,p_l=0}^{N-1} \mathop{
\sum_{m_1,n_1,\ldots,m_l,n_l=-\infty}}_{(\rref
{einschr1a}),\mbox{ }(\rref{einschr1b})}^\infty
\psi_{m_1}(u_{j_1})\cdots\psi_{n_l}(u_{j_l})
\exp \bigl(-\i \lambda_{k_1}(p_1-p_l)\bigr) \\
&& \hspace*{127.4pt}\quad{}\times
\prod
_{i=1}^{l-1} \exp\bigl(-\i
\lambda_{k_{i+1}}(p_{i+1}-p_i)\bigr)\\
&& \hspace*{127.4pt}\hspace*{22.4pt}\quad{}\times
\exp\bigl(-\i \lambda_{k_1}(m_l-n_1+t_{j_1}-t_{j_l})
\bigr) \\
&& \hspace*{127.4pt}\quad{}\times\prod_{i=1}^{l-1}\exp\bigl(-\i
\lambda_{k_{i+1}}(m_i-n_{i+1}+t_{j_{i+1}}-t_{j_i})
\bigr).
\end{eqnarray*}
We rename the $m_i, n_i$ ($m_i$ is replaced by $n_i$ and $n_i$ is
replaced with $m_{i-1}$ where we identify $l+1$ with $1$ and $0$ with
$l$). Then (\ref{einschr1a}) and (\ref{einschr1b}) become
%
\begin{eqnarray}
\label{einschr2a} 0 & \leq & p_i+m_i-n_i+t_{j_i}-t_{j_{i+1}}
\leq N-1 \qquad\mbox{for } i \in\{1,\ldots, l-1\},
\\
\label{einschr2b} 0 & \leq & p_{l}+m_{l}-n_{l}+t_{j_{l}}-t_{j_1}
\leq N-1
\end{eqnarray}
and after a factorisation in the arguments of the exponentials we
obtain that $V_T(\tilde\nu)$ is equal to
\begin{eqnarray*}
&&\frac{1}{2^lT^{l/2}}\sum_{j_1,\ldots,j_l=1}^M\sum
_{k_1,\ldots,k_l=-\lfloor
({N-1})/{2} \rfloor}^{{N}/{2}} \phi(u_{j_1},
\lambda_{k_1})\cdots \phi(u_{j_l},\lambda_{k_l})
\frac{1}{(2\uppi N)^l}\\
&&\quad\hspace*{0pt}{}\times\sum_{p_1,p_2,\ldots,p_l=0}^{N-1} \mathop{
\sum_{m_1,n_1,\ldots,m_l,n_l=-\infty}}_{ (\rref
{einschr2a}),\mbox{ }(\rref{einschr2b})}^\infty
\psi_{m_1}(u_{j_2})\cdots\psi_{n_l}(u_{j_l})\\
&&\quad\hspace*{127pt}{}\times\prod_{i=1}^{l-1}\exp \bigl(-\i (
\lambda_{k_i}-\lambda_{k_{i+1}})p_i\bigr) \exp
\bigl(-\i (\lambda_{k_{l}}-\lambda_{k_{1}})p_l\bigr)
\\
&&\hspace*{160pt}{} \times\exp\bigl(-\i \lambda_{k_1}(n_l-m_l+t_{j_1}-t_{j_l})
\bigr)\\
&&\quad\hspace*{127pt}{} \times
\prod_{i=1}^{l-1}\exp \bigl(-\i
\lambda_{k_{i+1}}(n_i-m_i+t_{j_{i+1}}-t_{j_{i}})
\bigr).
\end{eqnarray*}
We see that one can divide the sum with respect to $p_i,m_i,n_i$ into a
product of two sums, namely one sum with respect to all $p_i,m_i,n_i$
with even $i$ and the same sum with odd $i$. Analogously, we divide
(\ref{einschr2a}) and (\ref{einschr2b}) into
%
\begin{equation}
\label{einschr3a} 0 \leq p_i+m_i-n_i+t_{j_i}-t_{j_{i+1}}
\leq N-1 \qquad\mbox{for } i \in\{1,3,5,\ldots,l-3,l-1\}
\end{equation}
and
%
\begin{eqnarray}
\label{einschr4a} 0 & \leq & p_i+m_i-n_i+t_{j_i}-t_{j_{i+1}}
\leq N-1 \qquad\mbox{for } i \in\{2,4,6,\ldots,l-4, l-2\},
\\
\label{einschr4b} 0 & \leq & p_{l}+m_{l}-n_{l}+t_{j_{l}}-t_{j_1}
\leq N-1.
\end{eqnarray}
After applying the Cauchy--Schwarz inequality we obtain that $V_T(\tilde
\nu)$ is bounded by
%
\begin{eqnarray}\label{csu}
&& \Biggl\{ \frac{1}{2^lT^{l/2}} \sum_{j_1,\ldots,j_l=1}^M
\sum_{k_1,\ldots,k_l=-\lfloor({N-1})/{2} \rfloor}^{{N}/{2}} \phi (u_{j_1},
\lambda_{k_1})^2\phi(u_{j_3},\lambda_{k_3})^2
\cdots\phi (u_{j_{l-1}},\lambda_{k_{l-1}})^2
\frac{1}{(2\uppi N)^{l}}
\nonumber\\
&&\quad{}\times \Biggl|\sum_{p_1=0}^{N-1} \exp\bigl(-\i (
\lambda_{k_1}-\lambda_{k_2})p_1\bigr)\sum
_{p_3=0}^{N-1} \exp\bigl(-\i (\lambda_{k_3}-
\lambda_{k_4})p_3\bigr) \cdots\\
&&\hspace*{13.3pt}\quad{}\times \sum
_{p_{l-1}=0}^{N-1} \exp\bigl(-\i (\lambda_{k_{l-1}}-
\lambda_{k_{l}})p_{l-1}\bigr)
\nonumber\\
&&\hspace*{13.3pt}\quad{}\times \mathop{\sum_{m_1,n_1,m_3,n_3,\ldots,m_{l-1},n_{l-1}=-\infty}}_{
(\rref{einschr3a})}^\infty
\psi_{m_1}(u_{j_2})\psi_{n_1}(u_{j_1})
\psi_{m_3}(u_{j_4})\psi_{n_3}(u_{j_3})
\cdots\nonumber\\
&&\hspace*{116pt}\hspace*{13.3pt}\quad{}\times
\psi_{m_{l-1}}(u_{j_l})\psi_{n_{l}}(u_{j_{l-1}})
\nonumber\\
&&\hspace*{13.3pt}\quad{}\times\prod_{a \in\{1,3,\ldots,l-1\}} \exp\bigl(-\i
\lambda_{k_{a+1}}(n_a-m_a+t_{j_{a+1}}-t_{j_a})
\bigr) \Biggr|^2 \Biggr\}^{1/2}
\nonumber\\
&&\quad{} \times \{ \mbox{the same term with even } p_i,m_i,n_i
\}^{1/2}.\nonumber
\end{eqnarray}
We only consider the first term in (\ref{csu}), which is equal to
%
\begin{eqnarray}\label{J_T}
J_T&:=&\frac{1}{2^l T^{l/2} } \sum
_{j_1,\ldots,j_l=1}^M\sum_{k_1,\ldots,k_l=-\lfloor({N-1})/{2} \rfloor}^{{N}/{2}}
\phi (u_{j_1},\lambda_{k_1})^2
\phi(u_{j_3},\lambda_{k_3})^2\cdots\phi
(u_{j_{l-1}},\lambda_{k_{l-1}})^2\hspace*{28pt}
\nonumber\\[-8pt]\\[-8pt]
&&\hspace*{134.3pt}{} \times\frac{1}{(2\uppi N)^{l}} \bigl|K_T(u_1,\ldots,u_l,
\lambda_{k_1},\ldots,\lambda_{k_l})\bigr|^2\hspace*{28pt}\nonumber
\end{eqnarray}
with $K_T(u_1,\ldots,u_l,\lambda_{k_1},\ldots,\lambda_{k_l})$ being defined
implicitly. We have
\begin{eqnarray*}
&& \frac{1}{(2\uppi N)^l}\sum_{k_2,k_4,\ldots,k_l=-\lfloor({N-1}){2}
\rfloor}^{{N}/{2}}\bigl|K_T(u_1,\ldots,u_l,
\lambda_{k_1},\ldots,\lambda_{k_l})\bigr|^2
\\
&&\quad=\frac{1}{(2\uppi N)^l} \sum_{p_1,p_3,\ldots,p_{l-1}=0}^{N-1} \sum
_{\tilde
p_1,\tilde p_3,\ldots,\tilde p_{l-1}=0}^{N-1} \mathop{\sum
_{m_1,n_1,m_3,n_3,\ldots,m_{l-1},n_{l-1}=-\infty}}_{ (\rref
{einschr3a})}^\infty \mathop{\sum
_{\tilde m_1,\tilde n_1,\tilde m_3,\tilde n_3,\ldots,\tilde
m_{l-1},\tilde n_{l-1}=-\infty}}_{\tilde
{(\rref{einschr3a})}}^\infty
\\
&&{}\qquad\exp\bigl(-\i \lambda_{k_1}(p_1-\tilde p_1)
\bigr)\exp\bigl(-\i \lambda_{k_3}(p_3-\tilde p_3)
\bigr) \cdots\exp\bigl(-\i \lambda_{k_{l-1}}(p_{l-1}-\tilde
p_{l-1})\bigr)\\
&&\hspace*{0pt}\qquad{}\times
\psi_{m_1}(u_{j_2})\psi_{n_1}(u_{j_1})
\cdots
\psi_{m_{l-1}}(u_{j_l})\psi_{n_{l-1}}(u_{j_{l-1}})
\psi_{\tilde
m_1}(u_{j_2})\psi_{\tilde n_1}(u_{j_1})
\cdots\\
&&\hspace*{0pt}\qquad{}\times\psi_{\tilde
m_{l-1}}(u_{j_l})\psi_{\tilde n_{l-1}}(u_{j_{l-1}})
\\
&&\qquad{}\times\sum_{k_2,k_4,\ldots,k_l=-
\lfloor({N-1})/{2} \rfloor}^{{N}/{2}} \exp\bigl(-\i
\lambda_{k_2}(\tilde p_1-p_1+n_1-m_1+
\tilde m_1-\tilde n_1)\bigr)\\
&&\hspace*{91pt}\qquad{}\times\exp \bigl(-\i
\lambda_{k_4}(\tilde p_3-p_3+n_3-m_3+
\tilde m_3-\tilde n_3)\bigr)\cdots
\\
&&\hspace*{91pt}\qquad{}\times\exp\bigl(-\i \lambda_{k_l}(\tilde p_{l-1}-p_{l-1}+n_{l-1}-m_{l-1}+
\tilde m_{l-1}-\tilde n_{l-1})\bigr)
\end{eqnarray*}
and because of the well-known identity
\[
\frac{1}{N} \sum_{k=- \lfloor({N-1})/{2} \rfloor}^{{N}/{2}} \exp
(-\i  \lambda_k t) = %
\cases{ 1, &\quad $t=lN$ with $l \in
\mathbb{Z}$,
\cr
0, &\quad else,}
\]
it follows that for every $i$ only one of the $p_i$ and $\tilde p_i$
can be chosen freely if the $m_i,n_i$ are fixed. Furthermore, we can
show with the same arguments as in the proof of (\ref{cov}) that
because of (\ref{einschr3a}) and (\ref{2.1a}) we only have to consider
the cases with $j_i = j_{i+1}$ for every odd $i$ and that all other
terms are of order $\O(1/N)$. This implies
\[
\frac{1}{(2\uppi N)^l}\sum_{k_2,k_4,\ldots,k_l=-\lfloor({N-1})/{2}
\rfloor}^{{N}/{2}}\bigl|K_T(u_1,\ldots,u_l,
\lambda_{k_1},\ldots,\lambda_{k_l})\bigr|^2 \leq
\frac{1}{(2\uppi)^l}\Biggl(\sum_{m=-\infty}^\infty|
\psi_m|\Biggr)^{2l}
\]
with $|\psi|:=\sup_u |\psi(u)|$, and since we only need to sum over
$j_i$ with odd $i$ in (\ref{J_T}), it follows
\[
J_T\leq\frac{1}{T^{l/2}(4\uppi)^l} \Biggl(\sum_{m=-\infty}^\infty|
\psi_m|\Biggr)^{2l} \Biggl(\sum
_{j=1}^M \sum_{k=1}^{{N}/{2}}
\phi (u_{j},\lambda_{k})^2
\Biggr)^{l/2} \bigl(1+\O(1/N)\bigr).
\]
We obtain the same upper bound for the second factor in (\ref{csu}) and
this implies
\begin{eqnarray*}
\cum_l\bigl(\sqrt{T} \hat D_T(\phi)\bigr) &\leq& \sum
_{\nu} \frac{1}{(4\uppi
)^l(2\uppi
)^{l/2}} \Biggl(\sum
_{m=-\infty}^\infty|\psi_m|\Biggr)^{2l}
\biggl( \int_0^1\int_0^{\uppi}
\phi^2(u,\lambda) \,\d \lambda \,\d u \biggr)^{l/2}
\\
&&{}\times\bigl(1+\O\bigl(N^2/T^2\bigr)+\O(1/N)\bigr)
\\
&\leq& (2l)! 2^l \frac{1}{(4\uppi)^l(2\uppi)^{l/2}} \Biggl(\sum
_{m=-\infty
}^\infty |\psi_m|
\Biggr)^{2l} \biggl( \int_0^1\int
_0^{\uppi} \phi^2(u,\lambda) \,\d \lambda \,\d u
\biggr)^{l/2}
\\
&&{}\times\bigl(1+\O\bigl(N^2/T^2\bigr)+\O(1/N)\bigr)
\\
&\leq&(2l)!\tilde C^l \rho_2(\phi)^l,
\end{eqnarray*}
where the last inequality follows because of $N/T \rightarrow0$ and
$1/N \rightarrow0$ and since $(2l)! 2^l$ is an upper bound for the
number of indecomposable partitions of (\ref{part1})
(see Dahlhaus~\cite{dahlhaus1988}).
\end{pf*}
\noqed\end{pf*}

\subsection{\texorpdfstring{Proof of Theorem \protect\ref{thm2}}{Proof of Theorem 3.3}} \label{sec52}
A consequence of assumption (\ref{2.1a}) and $\int_0^1f(u,\lambda) \,\d u
>0$ for all $\lambda\in[-\uppi,\uppi]$ together with Lemma 2.1 of
Kreiss, Paparoditis and Politis \cite{kreisspappol2011} is that
%
\begin{equation}
\label{arbound} \sum_{j=1}^\infty j
|a_j|<\infty
\end{equation}
holds, and Lemma 2.3 in Kreiss, Paparoditis and Politis
\cite{kreisspappol2011} implies that there exists a $p_0 \in\N$ such
that for all $p \geq p_0$ the $\operatorname{AR}(p)$ process defined through
(\ref{arpstat}) has an $\operatorname{MA}(\infty)$ representation
%
\begin{equation}
\label{arpbest} Y_t^{\mathrm{AR}}(p)=\sum
_{l=0}^\infty\psi_l^{\mathrm{AR}}(p)
Z_{t-l}^{\mathrm{AR}}(p).
\end{equation}
Furthermore, (\ref{arestbound}) together with Lemma 2.3 in
Kreiss, Paparoditis and Politis \cite{kreisspappol2011} imply that there exist a $p_0' \in\N$, such that
for all $p \geq p_0'$ the fitted $\operatorname{AR}(p)$ process has an $\operatorname{MA}(\infty)$
representation
%
\begin{equation}
\label{bprocess} X_{t,T}^*=\sum_{l=0}^\infty
\hat\psi_l^{\mathrm{AR}}(p,T) Z_{t-l}^*,
\end{equation}
and we assume without loss of generality that $T$ and $p(T)$ are
sufficiently large to ensure the existence of such a representation.

Recall the proof of Theorem \ref{thm1}. In case the process is
stationary, all the terms of order
$\O(N^2/T^2)$ and $\O(1/T)$ vanish, as they are due to certain
approximation errors which do not appear for $\psi_{t,T,l}=\psi_l(u)=\psi_l$. For a fixed $p$ and $T$, the process of interest (\ref
{bprocess}) is now indeed a stationary one and therefore the proof of
Theorem \ref{thm2} works in the same way as the previous one, if the
remaining terms (which are the ones of order
$\O_P(1/N)$) are a $\o_P(T^{-1/2})$ for the bootstrap process as well.
Even more precisely, we only need the terms of order
$\O_P(1/N)$ to be a $\o_P(T^{-1/2})$ in the calculation of the
expectation, while it would suffice that they are a $\o_P(1)$ in the
calculation of the higher order cumulants. A detailed look at the proof
of Theorem \ref{thm1} reveals that these terms are up to a constant of
the form
\[
\frac{(\sum_{m=0}^\infty|\psi_m|)^{q_1}(\sum_{l=0}^\infty l |\psi_l|)^{q_2}}{N}
\]
with $q_1,q_2 \in\N$. For example, if the process is stationary we
obtain from (\ref{boundgrN}) an upper bound for $|\E(\hat
D_T(u,\lambda
))|$ via
\[
C\frac{\sum_{m=0}^\infty|\psi_m|\sum_{l=0}^\infty l |\psi_l| }{N}=\O(1/N)
\]
for some $C \in\R$, so an upper bound for the expectation of the
bootstrap analogue $\hat D_T^*(u,\lambda)$ of $\hat D_T(u,\lambda)$ is
given by
\[
C \frac{\sum_{m=0}^\infty|\hat\psi_m^{\mathrm{AR}}(p,T)| \sum_{l=0}^\infty l
|\hat\psi_l^{\mathrm{AR}}(p,T)|}{N}.
\]
Therefore, it needs to be shown that
\[
\sqrt{T} \frac{\sum_{m=0}^\infty|\hat\psi_m^{\mathrm{AR}}(p,T)| \sum_{l=0}^\infty l
|\hat\psi_l^{\mathrm{AR}}(p,T)|}{N}=\o_P(1)
\]
holds to obtain
\[
\sqrt{T}\E\bigl(\hat D_T^*(u,\lambda)\bigr)=\o_P(1).
\]
Because of (\ref{arestbound}), we can use the following bound from the
proof of Theorem 3.1 in Berg, Paparoditis and Politis
\cite{bergpappolitis2010} for the difference between
$\hat\psi_l^{\mathrm{AR}}(p,T)$ and $\psi_l^{\mathrm{AR}}(p)$ which is uniform in $p(T)$
and in $l \in\N$:
%
\begin{equation}
\label{boundappr1} \bigl|\hat\psi_l^{\mathrm{AR}}(p,T) -
\psi_l^{\mathrm{AR}}(p)\bigr| \leq p(1+1/p)^{-l}
\O_P(\sqrt {\log{T}/T}).
\end{equation}
With (\ref{boundappr1}), we obtain
\[
\sum_{l=0}^\infty\bigl|\hat\psi_l^{\mathrm{AR}}(p,T)
- \psi_l^{\mathrm{AR}}(p)\bigr|=\O_P\bigl(p_{\mathrm{max}}^2(T)
\sqrt{\log{T}/T}\bigr)
\]
and
\[
\sum_{l=0}^\infty l \bigl|\hat
\psi_l^{\mathrm{AR}}(p,T) - \psi_l^{\mathrm{AR}}(p)\bigr|
=\O_P\bigl(p_{\mathrm{max}}^3(T) \sqrt{\log{T}/T}\bigr)
\]
using properties of the geometric series, which yields
\[
\sum_{l=0}^\infty\bigl|\hat\psi_l^{\mathrm{AR}}(p,T)
\bigr| \leq \O_P\bigl(p_{\mathrm{max}}^2(T) \sqrt {\log{T}/T}
\bigr) + \sum_{l=0}^\infty\bigl|
\psi_l^{\mathrm{AR}}(p)\bigr|
\]
and
\[
\sum_{l=0}^\infty l \bigl|\hat
\psi_l^{\mathrm{AR}}(p,T) \bigr| \leq \O_P
\bigl(p_{\mathrm{max}}^3(T) \sqrt {\log{T}/T}\bigr) + \sum
_{l=0}^\infty l \bigl|\psi_l^{\mathrm{AR}}(p)\bigr|.
\]
Lemma 2.4 of Kreiss, Paparoditis and Politis \cite{kreisspappol2011} now implies that
%
\begin{equation}
\label{kineq} \sum_{l=1}^\infty(1+l)\bigl|
\psi_l^{\mathrm{AR}}(p)-\psi_l\bigr| \leq\tilde C \sum
_{l=p+1}^\infty(1+l) |a_l|
\end{equation}
for another constant $\tilde C \in\R$, where the $a_l$ are the
coefficients of the $\operatorname{AR}(\infty)$ representation in (\ref{statproc}).
Note that we implicitly assumed in (\ref{kineq}) that the $\psi_l$ are
the coefficients of the Wold representation of the process $Y_t$
defined in (\ref{statproc}), since this particular bound only holds for
this special $\mathrm{MA}$ representation. However, since the proof of Theorem
\ref{thm1} does not depend at all on the kind of $\mathrm{MA}$ representation,
we can assume without loss of generality that the $\psi_l$ are the
coefficients of the Wold representation, and then (\ref{kineq})
together with (\ref{2.1a}) and (\ref{arbound}) yields
\[
\sum_{l=0}^\infty l \bigl|\psi_l^{\mathrm{AR}}(p)\bigr|
\leq\bar C
\]
for $\bar C \in\R$. Therefore, we obtain with (\ref{asspNM})
\[
\Biggl( \sum_{m=0}^\infty\bigl|\hat
\psi_m^{\mathrm{AR}}(p,T) \bigr| \Biggr)^{p_1} \Biggl( \sum
_{l=0}^\infty l \bigl|\hat\psi_l^{\mathrm{AR}}(p,T)
\bigr| \Biggr)^{p_2}=\O_P(1)
\]
for $p_1, p_2 \in\N$, which yields the assertion.
\end{appendix}

\section*{Acknowledgements}

This work has been supported in part by the Collaborative Research
Center ``Statistical modeling of nonlinear dynamic processes'' (SFB
823, Teilprojekt A1, C1) of the German Research Foundation (DFG). The
authors would like to thank two referees and an Associate Editor for
their constructive comments on an earlier version of this manuscript.
We are also grateful to Martina Stein who typed parts of this
manuscript with considerable technical expertise.



\printhistory

\end{document}